\documentclass[11pt, a4paper]{article}
\usepackage{amsmath}
\usepackage{graphics,graphicx}
\usepackage{verbatim}
\usepackage{amssymb}
\usepackage{latexsym}
\usepackage{amsmath, amsthm, amssymb}

\newtheorem{corollary}{Corollary}[section]

\newtheorem{lemma}[corollary]{Lemma}
\newtheorem{proposition}[corollary]{Proposition}
\newtheorem{remark}[corollary]{Remark}
\newtheorem{theorem}[corollary]{Theorem}
\newcommand{\mylabel}[1]{\label{#1}
            \ifx\undefined\stillediting
            \else \fbox{$#1$}\fi }
\newcommand{\BE}{\begin{equation}}
\newcommand{\BEQ}[1]{\BE\mylabel{#1}}
\newcommand{\EEQ}{\end{equation}}
\newcommand{\rfb}[1]{\mbox{\rm
   (\ref{#1})}\ifx\undefined\stillediting\else:\fbox{$#1$}\fi}

\newfont{\Blackboard}{msbm10 scaled 1200}
\newcommand{\bl}[1]{\mbox{\Blackboard #1}}
\newfont{\roma}{cmr10 scaled 1200}

\def\CC{\rm \hbox{C\kern-.56em\raise.4ex
         \hbox{$\scriptscriptstyle |$}\kern+0.5 em }}

\catcode`\é=\active \def é{\' e}

\newcommand{\nline}  {{\bl N}}
\newcommand{\rline}  {{\bl R}}

\def\b{\beta}

\def\s{\sigma}

\def\cA{{\cal A}}

\def\cD{{\cal D}}

\def\cH{{\cal H}}

\newcommand{\mm}    {{\hbox{\hskip 0.5pt}}}

\newcommand{\bluff} {{\hbox{\raise 15pt \hbox{\mm}}}}

\newcommand{\FORALL} {{\hbox{$\hskip 11mm \forall \;$}}}

\newcommand{\rarrow} {{\,\rightarrow\,}}

\makeatletter
\def\section{\@startsection {section}{1}{\z@}{-3.5ex plus -1ex minus
    -.2ex}{2.3ex plus .2ex}{\large\bf}}
\makeatother
\def\be{\begin{equation}}
\def\ee{\end{equation}}

\def\ds{\displaystyle}
\begin{document}
\thispagestyle{empty}
\title{\bf Spectral analysis and stabilization of a chain of serially connected Euler-Bernoulli beams and strings}
\author{
Ka\"{\i}s Ammari
\thanks{D\'epartement de Math\'ematiques, Facult\'e des Sciences de Monastir, 5019 Monastir, Tunisie,
  e-mail: kais.ammari@fsm.rnu.tn,} 
 , \, Denis Mercier $^\dag$, \, Virginie R\'egnier \thanks{Univ Lille Nord de France, F-59000 Lille, France,
 UVHC, LAMAV, FR CNRS 2956, F-59313 Valenciennes, France,
email: denis.mercier@univ-valenciennes.fr, virginie.regnier@univ-valenciennes.fr}\\
 and 
 \\
 Julie Valein \thanks{Institut Elie Cartan Nancy (IECN), Nancy-Universit\'e \& INRIA (Project-Team CORIDA), 
F-54506  Vandoeuvre-l\`es-Nancy Cedex France,
email: julie.valein@iecn.u-nancy.fr}}
\date{}
\maketitle
%
%
\begin{quotation}
{\bf Abstract.} {\small We consider $N$ Euler-Bernoulli beams and $N$ strings alternatively connected to one another and forming a particular 
network which is a chain beginning with a string. 
We study two stabilization problems on the same network and the spectrum of the corresponding conservative system: the characteristic equation as well as its asymptotic behavior are given. 
We prove that the energy of the solutions of the first dissipative system tends to zero when the time tends to infinity under some irrationality assumptions of the length of the strings and beams. 
On another hand we prove a polynomial decay result
of the energy of the second system, independently of the length of the strings and beams, for all regular initial data. Our technique is based on a
frequency domain method and combines a contradiction argument with the multiplier technique to carry out a special analysis for the
resolvent.}
\end{quotation}
\section{Introduction} \label{secintro}

\setcounter{equation}{0}
We consider the evolution problems $(P_1)$ and $(P_2)$ described by the following systems of $2N$ equations : 
\begin{equation*}
\leqno(P_1) 
\left \{
\begin{array}{l}
(\partial_t^2 u_{2j-1}-\partial_x^2u_{2j-1})(t,x)=0,\, x\in(0,l_{2j-1}),\, t\in(0,\infty),\, j = 1,...,N,\\ 
(\partial_t^2 u_{2j} +\partial_x^4 u_{2j})(t,x)=0,\, x \in (0,l_{2j}),\, t\in(0,\infty),\,  j = 1,...,N,\\
u_1(t,0)=0,\ u_{2N}(t,l_{2N})=0,\, t\in(0,\infty),\\
\partial_x^2 u_{2j}(t,0)= \partial_x^2 u_{2j}(t,l_{2j})= 0,\, t\in(0,\infty),\, j = 1,...,N,\\
u_j(t,l_j)=u_{j+1}(t,0),\, t\in(0,\infty),\, j = 1,...,2N-1,\\
\partial_x^3 u_{2j}(t,0)+\partial_x u_{2j-1}(t,l_{2j-1})= {\bf -\partial_t u_{2j-1}(t,l_{2j-1})},\, t\in(0,\infty),\, j = 1,...,N,\\
\partial_x^3 u_{2j}(t,l_{2j})+\partial_x u_{2j+1}(t,0)= {\bf \partial_t u_{2j}(t,l_{2j})},\, t\in(0,\infty),\, j = 1,...,N,\\
u_j(0,x)=u_j^0(x),\ \partial_t u_j(0,x)=u_j^1(x), \,x\in(0,l_j),\,   j=1,...,2N,
\end{array}
\right.\\
\end{equation*}
and
\begin{equation*}
\leqno(P_2) 
\left \{
\begin{array}{l}
(\partial_t^2 u_{2j-1}-\partial_x^2u_{2j-1})(t,x)=0,\, x\in(0,l_{2j-1}),\, t\in(0,\infty),\, j = 1,...,N,\\ 
(\partial_t^2 u_{2j} +\partial_x^4 u_{2j})(t,x)=0,\, x \in (0,l_{2j}),\, t\in(0,\infty),\, j = 1,...,N,\\
u_1(t,0)=0,\ u_{2N}(t,l_{2N})=0,\ \partial_x^2 u_{2N}(t,l_{2N})=0, \, t\in(0,\infty),\\
\partial_x^2 u_{2j}(t,0)=  {\bf \partial^2_{tx} u_{2j}(t,0)},\, t\in(0,\infty),\, j = 1,...,N,\\
\partial_x^2 u_{2j}(t,l_{2j})=  {\bf - \partial^2_{tx} u_{2j}(t,l_{2j})},\, t\in(0,\infty),\, j = 1,...,N-1,\\
u_j(t,l_j)=u_{j+1}(t,0),\, t\in(0,\infty),\, j = 1,...,2N-1,\\
\partial_x^3 u_{2j}(t,0)+\partial_x u_{2j-1}(t,l_{2j-1})={\bf - \, \partial_t u_{2j-1}(t,l_{2j-1})},\, t\in(0,\infty),\, j = 1,...,N,\\
\partial_x^3 u_{2j}(t,l_{2j})+\partial_x u_{2j+1}(t,0)={\bf  \partial_t u_{2j+1}(t,0)},\ t\in(0,\infty),\, j = 1,...,N-1,\\
u_j(0,x)=u_j^0(x),\ \partial_t u_j(0,x)=u_j^1(x), \,x\in (0,l_j),\,   j=1,...,2N,
\end{array}
\right.\\
\end{equation*}
where $l_j > 0, \, \forall \, j=1,...,2N$.

Models of the transient behavior of some or all of the state variables describing the motion of flexible structures have been of great interest
in recent years, for details about physical motivation for the models, see \cite{banks2}, \cite{dagerzuazua}, \cite{lagnese} and the references therein. Mathematical
analysis of transmission partial differential equations is detailed in \cite{lagnese}.


 Let us first introduce some notation and definitions which will be used throughout the rest of the paper, in particular some which are linked to the notion of $C^{\nu }$- networks,
$\nu \in \nline$ (as introduced in \cite{jvb} and recalled in \cite{merreg1}). \\
Let $\Gamma$ be a connected topological graph embedded in $\rline^2$, with $2N$ edges ($N \in \nline^{*}$).  
Let $K=\{k_{j}\, :\, 1\leq j\leq 2N\}$ be the set of the edges of $\Gamma$. Each edge $k_{j}$ is a Jordan curve in $\rline^{2}$ and is assumed to be parametrized by its arc length $x_{j}$ such that
the parametrization $\pi _{j}\, :\, [0,l_{j}]\rightarrow k_{j}\, :\, x_{j}\mapsto \pi _{j}(x_{j})$ is $\nu$-times differentiable, i.e. $\pi _{j}\in C^{\nu }([0,l_{j}],\rline^{2})$ for all $1\leq j\leq 2N$. The length of the edge $k_j$ is $l_j>0$. 
The $C^{\nu}$- network $G$  associated with $\Gamma$ is then defined as the union $$G=\bigcup _{j=1}^{2N}k_{j}.$$

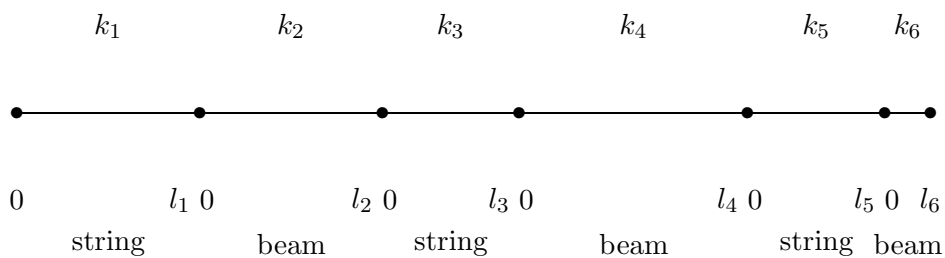
\begin{figure}[ht]
\begin{center}
\setlength{\unitlength}{0.240900pt}
\ifx\plotpoint\undefined\newsavebox{\plotpoint}\fi
\begin{picture}(1500,900)(0,0)
\sbox{\plotpoint}{\rule[-0.200pt]{0.400pt}{0.400pt}}%
\put(30,313){\makebox(0,0){$0$}}
\put(172,586){\makebox(0,0){$k_1$}}
\put(172,245){\makebox(0,0){string}}
\put(300,313){\makebox(0,0)[r]{$l_1$}}
\put(314,313){\makebox(0,0)[l]{$0$}}
\put(456,586){\makebox(0,0){$k_2$}}
\put(456,245){\makebox(0,0){beam}}
\put(583,313){\makebox(0,0)[r]{$l_2$}}
\put(598,313){\makebox(0,0)[l]{$0$}}
\put(704,586){\makebox(0,0){$k_3$}}
\put(704,245){\makebox(0,0){string}}
\put(796,313){\makebox(0,0)[r]{$l_3$}}
\put(810,313){\makebox(0,0)[l]{$0$}}
\put(988,586){\makebox(0,0){$k_4$}}
\put(988,245){\makebox(0,0){beam}}
\put(1151,313){\makebox(0,0)[r]{$l_4$}}
\put(1165,313){\makebox(0,0)[l]{$0$}}
\put(1272,586){\makebox(0,0){$k_5$}}
\put(1272,245){\makebox(0,0){string}}
\put(1364,313){\makebox(0,0)[r]{$l_5$}}
\put(1378,313){\makebox(0,0)[l]{$0$}}
\put(1414,586){\makebox(0,0){$k_6$}}
\put(1414,245){\makebox(0,0){beam}}
\put(1449,313){\makebox(0,0){$l_6$}}
\put(30,449){\usebox{\plotpoint}}
\put(30,449){\makebox(0,0){$\bullet$}}
\put(314,449){\makebox(0,0){$\bullet$}}
\put(598,449){\makebox(0,0){$\bullet$}}
\put(810,449){\makebox(0,0){$\bullet$}}
\put(1165,449){\makebox(0,0){$\bullet$}}
\put(1378,449){\makebox(0,0){$\bullet$}}
\put(1449,449){\makebox(0,0){$\bullet$}}
\put(30.0,449.0){\rule[-0.200pt]{341.837pt}{0.400pt}}
\end{picture}
\caption{A chain with $2N=6$ edges}
\end{center}
\label{fig}
\end{figure}

 We study two feedback stabilization problems for a string-beam network, see \cite{ammari1}-\cite{amjellk}, \cite{lagnese} and \cite{xumasto}-\cite{zhangxumasto}. In the following, only chains will be considered as mathematically described in Section 5 of \cite{merreg2}. See also \cite{merreg3} and Figure \ref{fig}. 


 Following Ammari/Jellouli/Mehrenberger (\cite{ammehjel}), we study a linear system modelling the vibrations of a chain of alternated Euler-Bernoulli beams and strings but with $N$ beams
and $N$ strings (instead of one string-one beam). For each edge $k_{j}$ (representing a string if $j$ is odd and a beam if $j$ is even), the scalar function $u_j(x,t)$ for $x \in G$ and $t > 0$
contains the information on the vertical displacement of the string if $j$ is odd and of the beam if $j$ is even ($1 \leq j \leq 2N$). 

Our aim is to study the spectrum of the conservative spatial operator 
which is defined in Section \ref{secspectrum} and to obtain stability results for $(P_1)$ and $(P_2)$.

We define the natural energy $E(t)$ of a solution $\underline{u} = (u_1,...,u_{2N})$ of $(P_1)$ or
$(P_2)$ by
\begin{multline}\label{energy1}
E(t)=\frac{1}{2} \ds \sum_{j=1}^{N} \left( \int_{0}^{l_{2j-1}} \left(|\partial_tu_{2j-1}(t,x)|^2+|\partial_xu_{2j-1}(t,x)|^2\right){\rm d}x\right.\\
\left.+ 
\int_{0}^{l_{2j}}\left(|\partial_tu_{2j}(t,x)|^2+|\partial_x^2u_{2j}(t,x)|^2\right){\rm d}x \right).
\end{multline}

We can easily check that every sufficiently smooth solution of $(P_1)$ satisfies the following dissipation law 
\begin{equation}\label{dissipae1}
E'(t) = - \ds \sum_{j=1}^{2N-1}\bigl|\partial_t u_{j}(t,l_{j})\bigr|^2\leq 0,
\end{equation}
and therefore, the energy is a nonincreasing function of the time variable $t$.

The first result concerns the well-posedness of the solutions of $(P_1)$ and the decay of the energy $E(t)$ of the solutions of $(P_1)$. 
We also study the spectrum of the corresponding conservative system. We give, in particular, the characteristic equation and the asymptotic behavior of the eigenvalues of the corresponding conservative system. We deduce that the generalized gap condition holds: if we denote by $(\lambda_n)_{n\in\mathbb{N}^*}$ the sequence of eigenvalues counted with their multiplicities, then 
\begin{equation}
\exists\gamma>0,\,\forall n\geq1,\,\lambda_{n+2N}-\lambda_n\geq\gamma.\label{generalgap}\end{equation}
Contrary to \cite{ammehjel}, it seems that the (simple) gap condition fails in general (for any $N\geq2$).
Therefore we do not succeed to obtain an observability inequality (and then to deduce stability results for $(P_1)$) directly by the study of the spectrum and the eigenvectors (see, for instance, \cite{valein:07}). In fact, the difficulties are to locate precisely the type of eigenvalues in the packets.

However, we prove that the energy $E(t)$ of the solutions of $(P_1)$ tends to zero when $t\rarrow +\infty$ in an appropriate energy space (described later), under some assumptions about the irrationality properties of the length of the strings and beams. For that, we use a result from \cite{arendt:88}.

\medskip

As we do not succeed to obtain the explicit decay rate to zero of the energy of the solutions of $(P_1)$, we change a little the system, by considering more dissipation conditions. That is why we introduce in problem $(P_2)$, in addition, the following dissipation conditions
$$\begin{array}{c}\partial_x^2 u_{2j}(t,0)=  {\bf \partial^2_{tx} u_{2j}(t,0)},\quad t\in(0,\infty),\quad j = 1,...,N,\\
\partial_x^2 u_{2j}(t,l_{2j})=  {\bf - \partial^2_{tx} u_{2j}(t,l_{2j})},\quad t\in(0,\infty), \quad j = 1,...,N-1.
\end{array}$$
In this case, we are able to prove more interesting stability results for system $(P_2)$ and to give the explicit decay rate of the energy of the solutions of $(P_2)$ in an appropriate space.

In the same manner as previously and with the same energy $E(t)$ (defined by (\ref{energy1})), every sufficiently smooth solution of $(P_2)$ satisfies the following dissipation law
\begin{equation}\label{energy}E'(t) = -\sum_{j=1}^{2N-1}\bigl|\partial_t u_{j}(t,l_j)\bigr|^2 
- \sum_{j=1}^{N-1} \bigl|\partial^2_{tx} u_{2j}(t,l_{2j})\bigr|^2 - \sum_{j=1}^{N}\bigl|\partial^2_{tx} u_{2j}(t,0)\bigr|^2\leq 0,
\end{equation}
and therefore, the energy is a nonincreasing function of the time variable $t$.

The main result of this paper then concerns the precise asymptotic behavior of the solutions of $(P_2)$.
As it was shown in \cite{ammehjel} in the case of one string and one beam connected together (i.e. $N=1$), we can not except to obtain an exponential decay rate of the solutions of $(P_2)$.
However we are able to prove that the decay rate to zero of the energy is $\ln^{4}(t)/t^2$, independently of the length of the strings and beams and by taking more regular initial data in an appropriate space.
Our technique is based on a frequency domain method from \cite{rao} and combines a contradiction argument with the multiplier technique to carry out a special analysis for the resolvent.

\medskip

This paper is organized as follows:
In Section \ref{wellposed}, we give the proper functional setting for systems $(P_1)$ and $(P_2)$ and prove that these two systems are well-posed.
In Section \ref{secspectrum}, we study the spectrum of the corresponding conservative system and we give the asymptotic behavior of the eigenvalues.
We then show that the energies of systems $(P_1)$ and $(P_2)$ tend to zero.
Finally, in Section \ref{end}, we study the stabilization result for $(P_2)$ by the frequency domain technique and give the explicit decay rate of the energy of the solutions of $(P_2)$.

\section{\label{wellposed}Well-posedness of the systems}

In order to study systems $(P_1)$ and $(P_2)$ we need a proper functional setting. 
We define the following space   
\begin{multline*}V= \bigg \{\underline{u}=(u_1,...,u_{2N}) \in \ds \prod_{j=1}^{N}\left( H^1(0,l_{2j-1}) \times H^2(0,l_{2j})\right),\\ 
u_j(l_j)=u_{j+1}(0), j=1, \ldots ,2N-1,\, u_1(0)=0,\, u_{2N}(l_{2N})=0 \bigg \},\end{multline*}
equipped with the sesquilinear form
\begin{equation}\label{ipV}
<\underline{u},\,\underline{\tilde{u}}>_{V}=\sum_{j=1}^{N}\left(\int_0^{l_{2j-1}}\partial_x u_{2j-1}(x)\partial_x  \overline{\tilde{u}_{2j-1}(x)}dx
+\int_0^{l_{2j}}\partial^2_x u_{2j}(x)\partial^2_x \overline{\tilde{u}_{2j}(x)}dx\right).
\end{equation}

Note the following lemma:
\begin{lemma}\label{lem0vp}
We have that $0$ is an eigenvalue associated to $(P_1)$ and $(P_2)$ of multiplicity $N-1$, i.e. there exists a subspace of $V$ of dimension $N-1$ such that any $\underline{\phi}$ in this subspace satisfies
\begin{equation*}
\leqno(EP_0) 
\left \{
\begin{array}{l}
\partial_x^2\phi_{2j-1}(x)=0,\, x\in(0,l_{2j-1}),\, j = 1,...,N,\\ 
\partial_x^4 \phi_{2j}(x)=0,\, x \in (0,l_{2j}),\,  j = 1,...,N,\\
\phi_1(0)=0,\ \phi_{2N}(l_{2N})=0,\\
\partial_x^2 \phi_{2j}(0)= \partial_x^2 \phi_{2j}(l_{2j})= 0,\,  j = 1,...,N,\\
\phi_j(l_j)=\phi_{j+1}(0),\, j = 1,...,2N-1,\\
\partial_x^3 \phi_{2j}(0)+\partial_x \phi_{2j-1}(l_{2j-1})= 0,\,  j = 1,...,N,\\
\partial_x^3 \phi_{2j}(l_{2j})+\partial_x \phi_{2j+1}(0)= 0,\, j = 1,...,N.
\end{array}
\right.\\
\end{equation*}
\end{lemma}

\begin{proof}
Let $\underline{\phi}$ be a non-trivial solution of $(EP_0)$. By the two first equations of $(EP_0)$, for $j\in\left\{1,\cdots,N\right\}$, $\phi_{2j-1}$ is a first order polynomial and $\phi_{2j}$ is a third order polynomial.
Moreover, with the fourth equation of $(EP_0)$, $\phi_{2j}$ also is a first order polynomial.
The two last equations of $(EP_0)$ become
$$\partial_x\phi_{2j-1}(0)=\partial_x\phi_{2j-1}(l_{2j-1})=0,\quad j=1,\cdots,N.$$
Consequently there exists $b_{2j-1}\in\mathbb{C}$ such that $\phi_{2j-1}=b_{2j-1}$ for $j\in\left\{1,\cdots,N\right\}$.
The third equation of $(EP_0)$ implies $b_1=0$.
Moreover we find, by the fifth equation of $(EP_0)$, that
$$\phi_{2j}(x)=\dfrac{b_{2j+1}-b_{2j-1}}{l_{2j}}\,x+b_{2j-1},\quad x\in(0,l_{2j}),\,j=1,\cdots,N,$$
where we set $b_{2N+1}=0$.

The function $\underline{\phi}$ defined above with $(b_3,b_5,\cdots,b_{2N-1})\in\mathbb{C}^{N-1}$ then satisfies $(EP_0)$, which finishes the proof.
\end{proof}

It is well-known that system $(P_1)$ may be rewritten as the first order evolution equation
\begin{equation} \left\{
\begin{array}{l}
U'=\mathcal{A}_1U,\\
U(0)=(\underline{u}^{0},\,\underline{u}^{1})=U_0,\end{array}\right.\label{pbfirstorder}\end{equation}
where $U$ is the vector $U=(\underbar{u},\,\partial_t \underbar{u})^t$ and the operator $\mathcal{A}_1:Y_1\rightarrow V\times\prod_{j=1}^{2N} L^2(0,l_{j})$ is defined by 
\[\mathcal{A}_1(\underline{u},\underline{v})^t:=(\underline{v}, (\partial_x^2u_{2j-1},-\partial_x^4 u_{2j})_{1 \leq j\leq N})^t,\] 
with
\begin{multline*}
Y_1:=\left\{(\underline{u},\,\underline{v})\in \prod_{j=1}^{N}\left( H^2(0,l_{2j-1}) \times H^4(0,l_{2j})\right)\times V \,:\right.\\
\left.\mbox {\textrm{satisfies }} \,(\ref {e2}) \; \mbox{\textrm{to}} \; (\ref {e5}) \; \mbox{\textrm{hereafter}}
\right\},\end{multline*}
\begin{equation}\label{e2}
\partial_x^2 u_{2N}(l_{2N}) = 0
\end{equation}
\begin{equation}\label{e3}
\partial_x^2 u_{2j}(0)= 0\quad j=1,...,N \quad \mbox {\textrm{and}} \quad 
\partial_x^2 u_{2j}(l_{2j})= 0, \quad j=1,...,N-1
\end{equation}
\begin{equation}\label{e4}
\partial_x^3 u_{2j}(0) + \partial_x u_{2j-1}(l_{2j-1})= -v_{2j-1}(l_{2j-1}), \quad j = 1,...,N 
\end{equation}
\begin{equation}\label{e5}
\partial_x^3 u_{2j}(t,l_{2j}) + \partial_x u_{2j+1}(0)=v_{2j}(l_{2j}), \quad j = 1,...,N-1.
\end{equation}

It is clear that $<\,.\,,\,.\,>_{ V}$ does not define a norm for $V$ but only a semi-norm since, for all $\underline{u}\in V$, we have 
$<{ \underline{u},\,\underline{u}}>_{V}=0$ if and only if $\underline{u}$ satisfies $(EP_0)$.
In order to get a Hilbert space we define by $E_0$, the eigenspace of $\mathcal{A}_1$ associated to the eigenvalue $0$, i.e.
$$E_0=\left\{(\underline{\phi},0)\in V\times \prod_{j=1}^{2N} L^2(0,l_{j})\,:\,\underline{\phi}\hbox{ satisfies } (EP_0)\right\},$$
and $P_{0,1}:V\times \prod_{j=1}^{2N} L^2(0,l_{j})\rightarrow E_0$ the projection onto $E_0$ defined by
$$P_{0,1}=\frac{1}{2\pi i}\oint_{\gamma}(\lambda I-\mathcal{A}_1)^{-1}d\lambda,$$
where $\gamma$ is a simple closed curve enclosing only the eigenvalue $0$ (see Theorem III-6.17 of \cite{Kato}).
Now let $\mathcal{H}_1$ the Hilbert space defined by
\begin{equation}\label{defH1}
V\times \prod_{j=1}^{2N} L^2(0,l_{j})=E_0\oplus\mathcal{H}_1,
\end{equation}
where $\mathcal{H}_1=(I-P_{0,1})(V\times \prod_{j=1}^{2N} L^2(0,l_{j}))$ and $E_0=P_{0,1}(V\times \prod_{j=1}^{2N} L^2(0,l_{j}))$.
Then $P_{0,1}$ is the projection onto $E_0$ parallel to $\mathcal{H}_1$. 
Note that, if $N=1$, $\mathcal{H}_1=V\times \prod_{j=1}^{2N} L^2(0,l_{j})$.

Then  $\mathcal{H}_1$ is a Hilbert space,
equipped with the usual inner product 
\begin{multline*}
\left\langle\left(\begin{array}{c}\underline{u}\\\underline{v}\end{array}\right),
\left(\begin{array}{c}\underline{\tilde{u}}\\ \underline{\tilde{v}}\end{array}\right)\right\rangle_{{\cal H}_1} = 
\sum_{j=1}^{N}\left(\int_{0}^{l_{2j-1}}\left(v_{2j-1}(x)\overline{\tilde{v}_{2j-1}(x)}
+\partial_x u_{2j-1}(x)\partial_x\overline{\tilde{u}_{2j-1}(x)}\right){\rm d}x\right.\\
\left.+\int_{0}^{l_{2j}}\left(v_{2j}(x)\overline{\tilde{v}_{2j}(x)}
+\partial_x^2u_{2j}(x)\partial_x^2\overline{\tilde{u}_{2j}(x)}\right){\rm d}x \right).\end{multline*}  

From now on we consider the operator $\mathcal{A}_1$ restricted to the space $\mathcal{H}_1\cap Y_1$ with value in $\mathcal{H}_1$, since $\mathcal{A}_1$ commutes with $P_{0,1}$.
By abuse of notation, this operator will be always denoted by $\mathcal{A}_1$ and 
 $\mathcal{D}(\mathcal{A}_1)$ will be its domain, i.e.
$$\mathcal{D}(\mathcal{A}_1):=\mathcal{H}_1\cap Y_1.$$

Therefore $$\mathcal{A}_1:\mathcal{D}(\mathcal{A}_1)\rightarrow\mathcal{H}_1.$$

 Moreover the norm on ${\cal D}({\cal A}_1)$ is defined by 
\be
\label{normDA}
\left| \left| (\underline{u},\underline{v}) \right| \right|_{{\cal D}({\cal A}_1)}^2 = \left| \left| {\cal A}_1 (\underline{u},\underline{v}) \right| \right|_{{\cal H}_1}^2 + \left| \left| 
(\underline{u},\underline{v})\right| \right|_{{\cal H}_1}^2.
\ee

\noindent Note that, with all these notation, problem $(P_1)$ is rewritten in an abstract way as: find $(\underline{u},\underline{v})^t \in {\cal D}({\cal A}_1)$ such that
$(\underline{u},\underline{v})^t_t = {\cal A}_1(\underline{u},\underline{v})^t$.

Now we can prove the well-posedness of system $(P_1)$ and that the solution of $(P_1)$ satisfies the dissipation law (\ref{dissipae1}).

\begin{proposition}\label{3exist1} 
(i) For an initial datum $U_{0}\in \mathcal{H}_1$, there exists a unique solution $U\in C([0,\,+\infty),\, \mathcal{H}_1)$
to  problem (\ref{pbfirstorder}). Moreover, if $U_{0}\in \mathcal{D}(\mathcal{A}_1)$, then
$$U\in C([0,\,+\infty),\, \mathcal{D}(\mathcal{A}_1))\cap C^{1}([0,\,+\infty),\, \mathcal{H}_1).$$ 

(ii) The solution $\underline{u}$ of $(P_1)$ with initial datum in $\mathcal{D}(\mathcal{A}_1)$ satisfies \rfb{dissipae1}.
Therefore the energy is decreasing.
\end{proposition}

\begin{proof}
(i) By Lumer-Phillips' theorem (see \cite{Pazy, tucsnakbook}), it suffices to show that
$\mathcal{A}_1$ is dissipative and  maximal.

We first prove that $\mathcal{A}_1$ is dissipative. Take $U=(\underline{u},\underline{v})^{t}\in \mathcal{D}(\mathcal{A}_1)$. 
Then
\begin{multline*}
\left\langle\mathcal{A}_1U,\, U \right\rangle_{\mathcal{H}_1}=\sum_{j=1}^{N}\left(\int_{0}^{l_{2j-1}}\left(\partial_x^2u_{2j-1}(x)\overline{v_{2j-1}(x)}
+\partial_x v_{2j-1}(x)\partial_x \overline{u_{2j-1}(x)}\right){\rm d}x\right.\\
\left.+\int_{0}^{l_{2j}}\left(-\partial_x^4 u_{2j}(x) \overline{v_{2j}(x)}
+\partial_x^2v_{2j}(x)\partial_x^2 \overline{u_{2j}(x)}\right){\rm d}x \right).
\end{multline*}
By integration by parts, we have
$$\Re\left(\left\langle\mathcal{A}_1U,\, U \right\rangle_{\mathcal{H}_1}\right)
=\Re\left(\sum_{j=1}^N\left[\partial_x u_{2j-1} \overline{v_{2j-1}}\right]_0^{l_{2j-1}}+\sum_{j=1}^N\left[-\partial_x^3 u_{2j} \overline{v_{2j}}\right]_0^{l_{2j}}
+\sum_{j=1}^N\left[\partial_x^2 u_{2j}\partial_x \overline{v_{2j}}\right]_0^{l_{2j}}\right).
$$
Moreover, we have
$$\sum_{j=1}^N\left[\partial_x^2 u_{2j}\partial_x \overline{v_{2j}}\right]_0^{l_{2j}}=0,
$$
by (\ref{e2}) and (\ref{e3}), and by the continuity of $\underline{v}$ at the interior nodes, we obtain
$$\begin{array}{l}
\displaystyle{\sum_{j=1}^N\left[\partial_x u_{2j-1} \overline{v_{2j-1}}\right]_0^{l_{2j-1}}
+\sum_{j=1}^N\left[-\partial_x^3 u_{2j} \overline{v_{2j}}\right]_0^{l_{2j}}}\\
=\displaystyle{\sum_{j=1}^N\left(\partial_x u_{2j-1}(l_{2j-1})+\partial_x^3 u_{2j}(0)\right)\overline{v_{2j-1}}(l_{2j-1})}
\displaystyle{-\sum_{j=1}^{N-1}\left(\partial_x u_{2j+1}(0)+\partial_x^3 u_{2j}(l_{2j})\right)\overline{v_{2j}}(l_{2j})}\\
\qquad\displaystyle{-\partial_xu_1(0)\overline{v_1}(0)-\partial_x^3u_{2N}(l_{2N})\overline{v_{2N}}(l_{2N})}\\
=\displaystyle{-\sum_{j=1}^N\left|v_{2j-1}(l_{2j-1})\right|^2-\sum_{j=1}^{N-1}\left|v_{2j}(l_{2j})\right|^2}
\end{array}
$$
by (\ref{e4}), (\ref{e5}) and since $\underline{v}\in V$.
Therefore
\begin{equation}\label{dissipativeness}
\Re\left(\left\langle\mathcal{A}_1U,\, U \right\rangle_{\mathcal{H}_1}\right)=-\sum_{j=1}^{2N-1}\left|v_{j}(l_{j})\right|^2\leq 0.
\end{equation}
This shows the dissipativeness of $\mathcal{A}_1$.

\medskip

Let us now prove that $\mathcal{A}_1$ is maximal, i.e. that
$\lambda I-\mathcal{A}_1$ is surjective for some $\lambda>0$.

Let $(\underline{f}, \underline{g})^{t}\in \mathcal{H}_1$. We look for $U=(\underline{u}, \underline{v})^{t}\in \mathcal{D}(\mathcal{A}_1)$ solution of 
\begin{equation}\label{eqmaxmon}
(\lambda I-\mathcal{A}_1)\left(\begin{array}{c}
\underline{u}\\\underline{v}\end{array}\right)=\left(\begin{array}{c}
\underline{f}\\\underline{g}\end{array}\right),\end{equation} 
or equivalently
\begin{equation} \left\{
\begin{array}{ll}
\lambda u_{j}-v_{j}=f_{j} & \forall j\in\{1,...,2N\},\\
\lambda v_{2j-1}-\partial^{2}_xu_{2j-1}=g_{2j-1} & \forall j\in\{1,...,N\},\\
\lambda v_{2j}+\partial^{4}_xu_{2j}=g_{2j}& \forall j\in\{1,...,N\}.\end{array}\right.\label{eqmaxmon2}\end{equation}

Suppose that we have found $\underline{u}$ with the appropriate regularity. 
Then for all  $j\in\{1,...,2N\},$ we have
\begin{equation} v_{j}:=\lambda u_{j}-f_{j}\in V.\label{maxmonv}\end{equation}
It remains to find $\underline{u}$. By (\ref{eqmaxmon2}) and (\ref{maxmonv}), $u_{j}$ must
satisfy, for all $j=1,...,N$,
$$ \lambda^{2}u_{2j-1}-\partial^{2}_xu_{2j-1}=g_{2j-1}+\lambda f_{2j-1},$$
and
$$\lambda^{2}u_{2j}+\partial^{4}_xu_{2j}=g_{2j}+\lambda f_{2j}.
$$ 
Multiplying these identities by a test function $\underline{\phi}$, integrating in space and using integration by
parts, we obtain
\begin{multline*}
\sum_{j=1}^{N}\int_0^{l_{2j-1}}\left(\lambda^2u_{2j-1}\overline{\phi_{2j-1}}+\partial_xu_{2j-1}\partial_x\overline{\phi_{2j-1}}\right)dx
-\sum_{j=1}^{N}\left[\partial_xu_{2j-1}\overline{\phi_{2j-1}}\right]_0^{l_{2j-1}}\\
+\sum_{j=1}^{N}\int_0^{l_{2j}}\left(\lambda^2u_{2j}\overline{\phi_{2j}}+\partial_x^2u_{2j}\partial_x^2\overline{\phi_{2j}}\right)dx
+\sum_{j=1}^{N}\left(\left[\partial_x^3u_{2j}\overline{\phi_{2j}}\right]_0^{l_{2j}}
-\left[\partial_x^2u_{2j}\partial_x\overline{\phi_{2j}}\right]_0^{l_{2j}}\right)\\
=\sum_{j=1}^{2N}\int_0^{l_j}\left(g_j+\lambda f_j\right)\overline{\phi_j}dx.
\end{multline*}
Since $(\underline{u},\underline{v})\in \mathcal{D}(\mathcal{A}_1)$ and $(\underline{u},\underline{v})$ satisfies (\ref{maxmonv}), we then have
\begin{multline}\label{maxmoneq1}
\sum_{j=1}^{N}\int_0^{l_{2j-1}}\left(\lambda^2u_{2j-1}\overline{\phi_{2j-1}}+\partial_xu_{2j-1}\partial_x\overline{\phi_{2j-1}}\right)dx
+\sum_{j=1}^{N}\int_0^{l_{2j}}\left(\lambda^2u_{2j}\overline{\phi_{2j}}+\partial_x^2u_{2j}\partial_x^2\overline{\phi_{2j}}\right)dx\\
+\sum_{j=1}^{2N-1}\lambda u_j(l_j)\overline{\phi_j}(l_j)=\sum_{j=1}^{2N}\int_0^{l_j}\left(g_j+\lambda f_j\right)\overline{\phi_j}dx
+\sum_{j=1}^{2N-1}f_j(l_j)\overline{\phi_j}(l_j).
\end{multline}
This problem has a unique solution $\underline{u}\in V$ by Lax-Milgram's lemma, because the left-hand side  of (\ref{maxmoneq1}) is coercive on $V$ equipped with the inner product defined by
\begin{multline*}
<\underline{u},\,\underline{\tilde{u}}>=\sum_{j=1}^{N}\int_0^{l_{2j-1}}(\partial_x u_{2j-1}(x)\partial_x  \overline{\tilde{u}_{2j-1}(x)}
+ u_{2j-1}(x) \overline{\tilde{u}_{2j-1}(x)})dx\\
+\sum_{j=1}^{N}\int_0^{l_{2j}}(\partial^2_x u_{2j}(x)\partial^2_x \overline{\tilde{u}_{2j}(x)}+ u_{2j}(x)\overline{\tilde{u}_{2j}(x)})dx,
\end{multline*}
and since $\lambda>0$.
If we consider $\underline{\phi}\in\prod_{j=1}^{2N}\mathcal{D}(0,\, l_{j})\subset V$, then $\underline{u}$
satisfies
$$\begin{array}{c}
\displaystyle{\lambda^{2}u_{2j-1}-\partial_x^{2}u_{2j-1}=g_{2j-1}+\lambda f_{2j-1} \quad\hbox{ in } \mathcal{D}'(0,\, l_{2j-1}),\quad  j=1,\cdots,N,}\\
\displaystyle{\lambda^{2}u_{2j}+\partial_x^{4}u_{2j}=g_{2j}+\lambda f_{2j} \quad \hbox{ in } \mathcal{D}'(0,\, l_{2j}),\quad  j=1,\cdots,N.}
\end{array}$$
This directly implies that  $\underline{u}\in\prod_{j=1}^{N}\left(H^{2}(0,\, l_{2j-1})\times H^{4}(0,\, l_{2j})\right)$ and then 
$\underline{u}\in V\cap \prod_{j=1}^{N}\left(H^{2}(0,\, l_{2j-1})\times H^{4}(0,\, l_{2j})\right)$. 
Coming back to (\ref{maxmoneq1}) and by integrating by parts, we find
$$\begin{array}{l}
\displaystyle{\sum_{j=1}^{N}\left(\partial_x^2u_{2j}(l_{2j})\partial_x\overline{\phi_{2j}}(l_{2j})
-\partial_x^2u_{2j}(0)\partial_x\overline{\phi_{2j}}(0)\right)}\\
\displaystyle{+\sum_{j=1}^{N}\left(\partial_xu_{2j-1}(l_{2j-1})+\partial_x^3u_{2j}(0)\right)\overline{\phi_{2j-1}}(l_{2j-1})}\\
\displaystyle{-\sum_{j=1}^{N-1}\left(\partial_xu_{2j+1}(0)+\partial_x^3u_{2j}(l_{2j})\right)\overline{\phi_{2j}}(l_{2j})
+\sum_{j=1}^{2N-1}\lambda u_j(l_j)\overline{\phi_j}(l_j)
=\sum_{j=1}^{2N-1}f_j(l_j)\overline{\phi_j}(l_j).}
\end{array}$$
Consequently, by taking particular test functions  $\underline{\phi}$, we obtain 
$$\begin{array}{c}
\partial_x^2 u_{2j}(l_{2j})=0 \quad \hbox{and}\quad \partial_x^2 u_{2j}(0)=0, \, j=1,\cdots,N,\\
\begin{array}{rcl}\partial_xu_{2j-1}(l_{2j-1})+\partial_x^3u_{2j}(0)&=&-\lambda u_{2j-1}(l_{2j-1})+f_{2j-1}(l_{2j-1})\\
	&=&-v_{2j-1}(l_{2j-1}),\, j=1,\cdots,N,\end{array}\\
\partial_xu_{2j+1}(0)+\partial_x^3u_{2j}(l_{2j})=\lambda u_{2j}(l_{2j})-f_{2j}(l_{2j})=v_{2j}(l_{2j}),\, j=1,\cdots,N-1.
\end{array}$$
It remains to show that  $(\underline{u},\underline{v})\in\mathcal{H}_1$.
Since $(\underline{f},\underline{g})\in\mathcal{H}_1=(I-P_{0,1})(V\times\prod_{j=1}^{2N}L^2(0,l_j))$, there exists $(\underline{\tilde f},\underline{\tilde g})\in V\times\prod_{j=1}^{2N}L^2(0,l_j)$ such that $(\underline{f},\underline{g})=(I-P_{0,1})(\underline{\tilde f},\underline{\tilde g}).$
Then 
$$(\underline{u},\underline{v})=(\lambda I-\mathcal{A}_1)^{-1}(\underline{f},\underline{g})=(\lambda I-\mathcal{A}_1)^{-1}(I-P_{0,1})(\underline{\tilde f},\underline{\tilde g})=(I-P_{0,1})(\lambda I-\mathcal{A}_1)^{-1}(\underline{\tilde f},\underline{\tilde g})\in\mathcal{H}_1,
$$ 
since the resolvent of $\mathcal{A}_1$ commutes with $P_{0,1}$ (see \cite{Kato}).

In summary we have found $(\underline{u},\underline{v})^{t}\in \mathcal{D}(\mathcal{A}_1)$ satisfying (\ref{eqmaxmon}), which finishes the proof of (i).

\medskip

(ii) To prove (ii), it suffices to derivate the energy (\ref{energy1})  for regular solutions and to use system $(P_1)$.
The calculations are analogous to those of the proof of the dissipativeness of $\mathcal{A}_1$ in (i), and then, are left to the reader.
\end{proof}

\bigskip

We see, in the same manner, that problem $(P_2)$ can be rewritten in an abstract way as: find $(\underline{u},\underline{v})^t \in {\cal D}({\cal A}_2)$ such that $(\underline{u},\underline{v})^t_t = {\cal A}_2(\underline{u},\underline{v})^t$, where
$\mathcal{A}_2:Y_2\rightarrow V\times\prod_{j=1}^{2N}L^2(0,l_j)$ for
\begin{multline*}
Y_2:=\left\{(\underline{u},\,\underline{v})\in \prod_{j=1}^{N}\left(H^2(0,l_{2j-1}) \times H^4(0,l_{2j})\right)\times V \,:\right.\\
\left.\mbox {\textrm{satisfies }} \,(\ref {e2}),\, (\ref {e4}),\, (\ref {e5})\, \mbox{\textrm{and}}\, (\ref{e6})\; \mbox{\textrm{hereafter}}
\right\},\end{multline*}
\begin{equation}\label{e6}
\partial_x^2 u_{2j}(0)= \partial_{x}v_{2j}(0),\, j=1,...,N  \mbox {\textrm{ and }} 
\partial_x^2 u_{2j}(l_{2j})= -\partial_{x}v_{2j}(l_{2j}), \, j=1,...,N-1,
\end{equation}
\[\mathcal{A}_2(\underline{u},\underline{v})^t:=(\underline{v}, (\partial_x^2u_{2j-1},-\partial_x^4 u_{2j})_{1 \leq j\leq N})^t.\]
Then we define the Hilbert space $\mathcal{H}_2$ by
$$V\times\prod_{j=1}^{2N}L^2(0,l_j)=E_0\oplus\mathcal{H}_2,\quad \mathcal{H}_2=(I-P_{0,2})(V\times\prod_{j=1}^{2N}L^2(0,l_j))$$
with $P_{0,2}:V\times\prod_{j=1}^{2N}L^2(0,l_j)\rightarrow E_0$ the projection onto $E_0$ defined by
$$P_{0,2}=\frac{1}{2i\pi}\oint_{\gamma}(\lambda I-\mathcal{A}_2)^{-1}d\lambda$$
(with $\gamma$ is a simple closed curve enclosing only the eigenvalue $0$), and
$$\mathcal{D}(\mathcal{A}_2):=\mathcal{H}_2\cap Y_2.
$$
Then $$\mathcal{A}_2:\mathcal{D}(\mathcal{A}_2)\rightarrow\mathcal{H}_2.$$
The following proposition holds:
\begin{proposition}\label{exist2}
(i) For an initial datum $U_{0}\in \mathcal{H}_2$, there exists a unique solution $U\in C([0,\,+\infty),\, \mathcal{H}_2)$ to  
$$ \left\{\begin{array}{l}
U'=\mathcal{A}_2U,\\
U(0)=(\underline{u}^{0},\,\underline{u}^{1})=U_0.\end{array}\right.
$$
Moreover, if $U_{0}\in \mathcal{D}(\mathcal{A}_2)$, then $$U\in C([0,\,+\infty),\, \mathcal{D}(\mathcal{A}_2))\cap C^{1}([0,\,+\infty),\, \mathcal{H}_2).$$ 

(ii) The solution $\underline{u}$ of $(P_2)$ with initial datum in $\mathcal{D}(\mathcal{A}_2)$ satisfies \rfb{energy}.
Therefore the energy is decreasing.
\end{proposition}

\begin{proof}
The proof of (i) and (ii) is the same as the proof of Proposition \ref{3exist1}, and therefore is left to the reader.
\end{proof}

%
%
\section{\label{secspectrum}Spectral analysis of a chain of serially connected Euler-Bernoulli beams and strings}

In this section, we study the spectral analysis of the corresponding conservative system.

Let $\underline{\Phi}$ be the solution of the conservative system derived from problems $(P_1)$ and $(P_2)$ given in the introduction, i.e. $\underline{\Phi}$ is the solution of the following system
\begin{equation*}
\leqno(P_c) 
\left \{
\begin{array}{l}
(\partial_t^2 \Phi_{2j-1}-\partial_x^2\Phi_{2j-1})(t,x)=0,\, x\in(0,l_{2j-1}),\, t\in(0,\infty),\, j = 1,...,N,\\ 
(\partial_t^2 \Phi_{2j} +\partial_x^4 \Phi_{2j})(t,x)=0,\, x \in (0,l_{2j}),\, t\in(0,\infty),\,  j = 1,...,N,\\
\Phi_1(t,0)=0,\ \Phi_{2N}(t,l_{2N})=0,\, t\in(0,\infty)\\
\partial_x^2 \Phi_{2j}(t,0)= \partial_x^2 \Phi_{2j}(t,l_{2j})= 0,\, t\in(0,\infty),\, j = 1,...,N,\\
\Phi_j(t,l_j)=\Phi_{j+1}(t,0), t\in(0,\infty),\, j = 1,...,2N-1,\\
\partial_x^3 \Phi_{2j}(t,0)+\partial_x \Phi_{2j-1}(t,l_{2j-1})= 0,\, t\in(0,\infty),\, j = 1,...,N,\\
\partial_x^3 \Phi_{2j}(t,l_{2j})+\partial_x \Phi_{2j+1}(t,0)= 0,\, t\in(0,\infty),\, j = 1,...,N,\\
\Phi_j(0,x)=u_j^0(x),\ \partial_t \Phi_j(0,x)=u_j^1(x), \, x\in(0,l_j),\,   j=1,...,2N,
\end{array}
\right.\\
\end{equation*}
where we have replaced the dissipative conditions (in bold in systems $(P_1)$ and $(P_2)$) by the conservative ones.

We can rewrite system $(P_c)$ in an abstract way as: find $(\underline{\Phi},\underline{\Psi})^t \in {\cal D}({\cal A}_c)$ such that $(\underline{\Phi},\underline{\Psi})^t_t = {\cal A}_c(\underline{\Phi},\underline{\Psi})^t$, where 
$\mathcal{A}_c:Y_c\rightarrow V\times\prod_{j=1}^{2N}L^2(0,l_j),$ for
\begin{multline*}
Y_c:=\left\{(\underline{\Phi},\,\underline{\Psi})\in\prod_{j=1}^{N}\left( H^2(0,l_{2j-1}) \times H^4(0,l_{2j})\right)\times V \,:\right.\\
\left.\mbox {\textrm{satisfies }} \,(\ref {e2}),\, (\ref {e3}),\, \mbox{\textrm{and}}\, (\ref{e7}),\, (\ref {e8})\; \mbox{\textrm{hereafter}}
\right\},\end{multline*}
\begin{equation}\label{e7}
\partial_x^3 \Phi_{2j}(0) + \partial_x \Phi_{2j-1}(l_{2j-1})= 0, \quad j = 1,...,N 
\end{equation}
\begin{equation}\label{e8}
\partial_x^3 \Phi_{2j}(t,l_{2j}) + \partial_x \Phi_{2j+1}(0)=0, \quad j = 1,...,N-1,
\end{equation}
and
\[\mathcal{A}_c(\underline{\Phi},\underline{\Psi})^t:=(\underline{\Psi}, (\partial_x^2\Phi_{2j-1},-\partial_x^4 \Phi_{2j})_{1 \leq j\leq N})^t.\]
Then we define the Hilbert space $\mathcal{H}_c$ by
\begin{equation}\label{defHc}
V\times\prod_{j=1}^{2N}L^2(0,l_j)=E_0\oplus\mathcal{H}_c,\quad \mathcal{H}_c=(I-P_{0,c})(V\times\prod_{j=1}^{2N}L^2(0,l_j))
\end{equation}
with $P_{0,c}:V\times\prod_{j=1}^{2N}L^2(0,l_j)\rightarrow E_0$ the projection onto $E_0$ defined by
with $$P_{0,c}=\frac{1}{2i\pi}\oint_{\gamma}(\lambda I-\mathcal{A}_c)^{-1}d\lambda$$
(with $\gamma$ is a simple closed curve enclosing only the eigenvalue $0$),
and
$$\mathcal{D}(\mathcal{A}_c):=\mathcal{H}_c\cap Y_c.
$$
Due to (\ref{defHc}), we set $V_c$ the Hilbert space defined by
$$\mathcal{H}_c=V_c\times\prod_{j=1}^{2N}L^2(0,l_j),$$
equipped with the inner product (\ref{ipV}).

Following Section \ref{wellposed}, 
it is clear that system $(P_c)$ is well-posed in the natural energy space.
If we suppose that $(\underline{u}^0,\underline{u}^1)\in \mathcal{H}_c=V_c\times\prod_{j=1}^{2N}L^2(0,l_j)$, then problem $(P_c)$ admits a unique solution
$$\underline{\Phi}\in C([0,T],V_c)\cap C^1([0,T],\prod_{j=1}^{2N}L^2(0,l_j)).
$$
This system is obviously conservative, i.e. its energy is constant.

\subsection{The characteristic equation}

Let $\underline{\phi}$ be a non-trivial solution of the eigenvalue problem $(EP)$ associated to the conservative problem $(P_c)$
and $\lambda^2$ be the corresponding eigenvalue. That is to say, $\underline{\phi}\in V_c$ satisfies the transmission and boundary conditions (\ref{p2})-(\ref{p6}) hereafter as well as 
$$(EP) \left \{ \begin{array}{ll}
\displaystyle \partial_x^2 \phi_{2j-1} = \lambda^{2} \phi_{2j-1} & \textrm{on}\, \textrm{ }(0,\textrm{l}_{2j-1}),\, \textrm{ }\forall \textrm{ j}\in \textrm{ }\{1,...,\textrm{N}\},\\
- \displaystyle \partial_x^4 \phi_{2j} = \lambda^{2} \phi_{2j} & \textrm{on}\, \textrm{ }(0,\textrm{l}_{2j}),\, \textrm{ }\forall \textrm{ j}\in \textrm{ }\{1,...,\textrm{N}\},\\
\phi_{2j-1} \in H^{2}(0,l_{2j-1}), & \forall j \in \{1,...,N\}, \; \phi_{2j}\in H^{4}(0,l_{2j}),\, \forall j \in \{1,...,N\}, \; 
\end{array} \right.$$
\be\label{p2}
\phi_1(0)=0, \; \phi_{2N}(l_{2N})=0,
\ee
\be\label{p3}
\partial_x^2 \phi_{2j}(0)= \partial_x^2 \phi_{2j}(l_{2j})= 0,\quad  j = 1,...,N
\ee
\be\label{p4}
\phi_j(l_j)=\phi_{j+1}(0),\quad j = 1,...,2N-1
\ee
\be\label{p5}
\partial_x^3 \phi_{2j}(0)+\partial_x \phi_{2j-1}(l_{2j-1})= 0,\quad j = 1,...,N
\ee
\be\label{p6}
\partial_x^3 \phi_{2j}(l_{2j})+\partial_x \phi_{2j+1}(0)= 0,\quad j = 1,...,N-1.
\ee

Note that this also means that $(\underline{\phi}, \lambda\underline{\phi})\in \mathcal{D}(\mathcal{A}_c)$ is an eigenvector of ${\cal A}_c$ associated to the eigenvalue $\lambda$.
By the definition of $\mathcal{A}_c$ and of its domain, $0$ is not an eigenvalue of $\mathcal{A}_c$. Moreover $0$ is not an eigenvalue of $\mathcal{A}_1$ and $\mathcal{A}_2$.

Define $z$ by $\lambda = i z^2$ where $z$ lies in $\rline^{+*}$ with $i^2 = - 1$. \\
Following Paulsen (\cite{paulsen}) and Mercier (\cite{mer}), we will rewrite this eigenvalue problem on a chain of $2N$ beams and 
strings using only square matrices of order $2$ in the following way: we define, for each $j\in\{1,...,N \}$, the vector functions $V_{2j-1}$ and $V_{2j}$ by 
$$V_{2j-1}(x)= \left ( \phi_{2j-1}(x),\, \dfrac{1}{z^2} \partial_x \phi_{2j-1}(x) \right)^{t},\,\forall x\in[0,l_{2j-1}],$$
$$V_{2j}(x)= \left ( \phi_{2j}(x),\, \dfrac{1}{z^3} \partial^3_x \phi_{2j}(x) \right)^{t} ,\,\forall x\in[0,l_{2j}].$$

\noindent Define the matrices $A_{j}$ by 
$$A_{2j-1}(z,l_{2j-1}):= \left(\begin{array}{cc}
c_{2j-1} & s_{2j-1} \\
-s_{2j-1} & c_{2j-1} 
\end{array}\right),$$
\begin{multline*}A_{2j}(z,l_{2j}):= \dfrac{1}{e^{2 l_{2j}z} - 2 e^{l_{2j}z} s_{2j} - 1}\\
 \cdot\left( \begin{array}{cc}
e^{2 l_{2j}z} \left( c_{2j} - s_{2j} \right) - c_{2j} - s_{2j} & 2 s_{2j} \left( 1 - e^{2 l_{2j}z} \right) \\
e^{2 l_{2j}z} c_{2j} - 2 e^{l_{2j}z} + c_{2j} & e^{2 l_{2j}z} \left( c_{2j} - s_{2j} \right) - c_{2j} - s_{2j} 
\end{array} \right),\end{multline*}
\noindent with $j\in\left\{1,\cdots,N\right\}$ and with the notation 
\begin{equation}\label{notcs}
\left\{ \begin{array}{cc}
c_{2j-1}= \cos(l_{2j-1} \cdot z^2), & s_{2j-1}= \sin(l_{2j-1} \cdot z^2) \\
c_{2j}= \cos(l_{2j} \cdot z), & s_{2j}= \sin(l_{2j} \cdot z).
\end{array} \right.\end{equation} 
 
\noindent The matrix $T$ is defined by:
$$T(z) := \left(\begin{array}{cc}
1 & 0 \\
0 & -\dfrac{1}{z} \end{array}\right).$$

\noindent To finish with, the matrix $M(z)$ is the square matrix of order $2$ given by  
\begin{equation}
M(z)=A_{2N}T A_{2N-1}...T^{-1} A_{2} T A_{1}.\label{eqM}
\end{equation}

\begin{lemma} (A few trivial but useful properties) \label{lemV}\\
With the notation introduced above, we have:
$$\begin{array}{lll}
V_{j}(l_{j})=A_{j}V_{j}(0),\, \forall j \in \{1,...,2N \}, \\
V_{2j}(0)=T V_{2j-1}(l_{2j-1}),\, \forall j \in \{1,...,N \}, \\
V_{2j+1}(0)= T^{-1} V_{2j}(l_{2j}),\, \forall j \in \{1,...,N-1 \}, \\
V_{2N}(l_{2N})=M(z)V_{1}(0).
\end{array}$$
\end{lemma}

\begin{proof} First, for $j$ odd and $j \in \{ 1, \ldots, 2N \}$, since $u_{j}$ satisfies the first equation of the eigenvalue problem $(EP)$, 
$u_{j}$ is a linear combination of the vectors of the fundamental basis
$$\left( \cos(z^2 \,.),\, \sin(z^2 \,.) \right).$$
\noindent The first equation of the lemma follows from that property after some calculations. \\
Now, for $j$ even and $j \in \{ 1, \ldots, 2N \}$, since $u_{j}$ satisfies the second equation of the eigenvalue problem $(EP)$, $u_{j}$ is a linear combination of the vectors of the fundamental basis
$$\left( \cos(z.),\, \sin(z\,.),\, e^{z\,.},\, e^{-z\,.} \right).$$ In this basis, if we consider the two following functions $d_1,d_2$ with coordinates 
$$d_1:=(-e^{l_jz}\sin(l_j z),e^{l_jz}\cos(l_jz)-1,0,-e^{l_jz}\sin(l_jz))$$
$$d_2:=(e^{l_j z}-e^{-l_j z},0,\cos(l_jz)-e^{-l_jz},e^{l_jz}-\cos(l_jz),$$ 
we can see that they are independent and satisfy (\ref{p2}). 
Consequently $u_j$ can be expressed as a linear combination of these two functions. Now, to find $A_j$, we proceed as follows: 
let $(\alpha,\beta)^t$ the coordinates of $u_j$ in the basis $(d_1,d_2)$. There exist two matrices $M_0, M_1$ such 
that $V_j(0)=M_0 (\alpha,\beta)^t$ and $V_j(l_j)=M_1 (\alpha,\beta)^t$, then $A_j$ is the matrix $M_1M_0^{-1}.$

Moreover the transmission conditions (\ref{p4}), (\ref{p5}) and (\ref{p6}) imply the second and third equations. \\
The fourth one is the logical consequence of the first three applied successively for $j=1$, $j=2$, etc... 
\end{proof}


\begin{theorem} (The characteristic equation for the eigenvalue problem corresponding to a chain of alternated beams and strings) \label{chareq} \\
The complex number $\lambda = i z^2$ ($z\in\mathbb{R}^{+*}$) is an eigenvalue of ${\cal {A}}_c$ if and only if $z$ satisfies the characteristic equation 
\begin{equation}
f(z)= m_{12}(z)=0,\label{car1}
\end{equation}
\noindent where $m_{12}(z)$ is the term on the first line and second column of the matrix $M(z)$.
\end{theorem}

\begin{proof} Let $\underline{\phi}$ be a non-trivial solution of the eigenvalue problem $(EP)$ and $\lambda^2$ be the corresponding eigenvalue, where $\lambda=iz^2$ ($z\in\mathbb{R}^{+*}$). 

\noindent Using the boundary conditions as well as $V_{2N}(l_{2N})=M(z)V_{1}(0)$, it follows: 
$$\left(\begin{array}{c}
0 \\ \dfrac{1}{z^3} \partial^3_x \phi_{2N}(l_{2N}) \end{array} \right)=M(z)\left(\begin{array}{c}
0 \\
\dfrac{1}{z^2} \partial_x \phi_{1}(0) 
\end{array}\right).$$
\noindent It is clear that the vector of the second part of the previous equality is non-trivial since $\underline{\phi}$ is a non-trivial
solution of problem $(EP)$. Hence the result. 
\end{proof}

\begin{proposition} (Asymptotic behavior of the characteristic equation) \label{asbehaviorspec}\\
Assume that the characteristic equation is given by Theorem \ref{chareq}. Then 
$$f(z)= z \left( f_{\infty}(z) + g(z) \right)$$
\noindent where
\begin{equation}
f_{\infty}(z) = s_1(z) \cdot c_2(z) \cdot s_3(z) \cdots s_{2N-1}(z) \cdot (c_{2N}(z) - s_{2N}(z)) 
\end{equation}
\noindent (with $c_j$, $s_j$ defined by (\ref{notcs})) and $g$ satisfies $\lim_{z \rightarrow +\infty} g(z)=0$. 
Thus, the asymptotic behavior of the spectrum $\sigma(\cA_c)$ corresponds to the roots of the asymptotic characteristic
equation 
\begin{equation}
f_{\infty}(z) = 0. \label{car2}
\end{equation}
\end{proposition}

\begin{proof}
In the following, the notation $o(h(\lambda))$ is used for a square matrix of order $2$ such that all its terms are dominated by
the function $\lambda \mapsto h(\lambda)$ asymptotically. For any $j \in \{ 1, \ldots, N \}$,
$$A_{2j}(z,l_{2j})= \dfrac{1}{e^{2 l_{2j}z} - 2 e^{l_{2j}z} s_{2j} - 1} \left[ e^{2 l_{2j}z} 
\left( \begin{array}{cc}
c_{2j} - s_{2j}  & 2 s_{2j} \\
c_{2j}           & c_{2j} - s_{2j} 
\end{array} \right) + o(2 l_{2j} z) \right].$$

Thus 
$$A_{2j}(z,l_{2j})= 
\left( \begin{array}{cc}
c_{2j} - s_{2j}  & 2 s_{2j} \\
c_{2j}           & c_{2j} - s_{2j} 
\end{array} \right) + o(1),$$
which leads, after some calculations, to:
$$T^{-1} A_{2j} T A_{2j-1} = 
\left( \begin{array}{cc}
(c_{2j} - s_{2j}) c_{2j-1}  & (c_{2j} - s_{2j}) s_{2j-1}   \\
-z c_{2j} c_{2j-1}          & -z c_{2j} s_{2j-1}           
\end{array} \right) + o(1). $$

Likewise
$$T^{-1} A_{2j+2} T A_{2j+1} = 
\left( \begin{array}{cc}
(c_{2j+2} - s_{2j+2}) c_{2j+1}  & (c_{2j+2} - s_{2j+2}) s_{2j+1}   \\
-z c_{2j+2} c_{2j+1}            & -z c_{2j+2} s_{2j+1}           
\end{array} \right) + o(1). $$

Thus 
\begin{multline*}
T^{-1} A_{2j+2} T A_{2j+1} T^{-1} A_{2j} T A_{2j-1} \\
= \left( \begin{array}{cc}
-z (c_{2j+2} - s _{2j+2}) s_{2j+1} c_{2j} c_{2j-1} & -z (c_{2j+2} - s_{2j+2}) s_{2j+1} c_{2j} s_{2j-1} \\
z^2 c_{2j+2} s_{2j+1} c_{2j} c_{2j-1}              &  z^2 c_{2j+2} s_{2j+1} c_{2j} s_{2j-1}                                  
\end{array} \right) + o(1).
\end{multline*}

The result follows by induction. \end{proof}

\begin{remark}
We can note that the eigenvalues $\lambda=iz^2$ of $(EP)$ have $2N$ families of asymptotic behavior: 
$$\begin{array}{c}
\ds{\left(i\,\frac{k\pi}{l_{2j-1}}\right)_{k\in\mathbb{N}^*},\quad j=1,\cdots,N,\quad \left(i\left(\frac{\pi+2k\pi}{2l_{2j}}\right)^2\right)_{k\in\mathbb{N}^*},\quad j=1,\cdots,N-1,}\\
\ds{\hbox{and } \left(i\left(\frac{\pi/4+k\pi}{l_{2N}}\right)^2\right)_{k\in\mathbb{N}^*}.}
\end{array}$$
It follows that the generalized gap condition (\ref{generalgap}) holds.
\end{remark}

\begin{proposition} (Geometric multiplicity of the eigenvalues) \label{mult}\\
If $\lambda\neq0$ is an eigenvalue of the operator ${\cal {A}}_c$ and $E_{\lambda}$ is the associated eigenspace, then the
dimension of $E_{\lambda}$ is one.
\end{proposition}

\begin{proof}
The eigenvectors $\underline{\phi}\in V_c$ associated to the eigenvalue $\lambda^2$ (cf. problem $(EP)$) are entirely determined by their 
values at the nodes of the network (i.e. where the beams and strings are connected to one another). Due to Lemma \ref{lemV}, they are
also determined by $V_1(0)= \left ( \phi_1(0),\, \dfrac{1}{z^2} \partial_x \phi_1(0) \right)^{t}$. Now $\phi_1(0) = 0$
(cf. condition ($\ref{p2}$)) and $\partial_x \phi_1(0)$ may take any value in $\rline^*$. Hence the result. \end{proof}

\subsection{Strong stability of $(P_1)$ and $(P_2)$}

We first prove the following lemma:
\begin{lemma}\label{lemmsum}
If there exist $i,\,j\in\left\{1,\cdots,N\right\}$ such that
 \begin{equation}\label{condstab1}
\frac{l_{2i-1}}{l_{2j-1}}\notin\mathbb{Q}\quad\hbox{or}\quad \frac{l_{2i}}{l_{2j}}\notin\mathbb{Q},
\end{equation}
or if there exist $i,\,j\in\left\{1,\cdots,N\right\}$ such that
 \begin{equation}\label{condstab2}
\frac{(l_{2i})^2}{l_{2j-1}}\neq\frac{p^2}{q}\pi,\quad\hbox{where}\quad p,\,q\in\mathbb{Z},
\end{equation}
then
\begin{equation}\label{sumneq0}
\sum_{j=1}^{2N-1}\left|\phi_j(l_j)\right|^2\neq 0,
\end{equation}
for all eigenvectors $\underline{\phi}\in V_c$ of $(EP)$. 
\end{lemma}

\begin{proof}
Let $\underline{\phi}\in V_c$ be an eigenvector of $(EP)$ associated to the eigenvalue $\lambda^2$, where $\lambda=iz^2$ ($z\in\mathbb{R}^{+*}$).
Assume that (\ref{sumneq0}) is false, i.e. that we have
\begin{equation}\label{sum=0}
\sum_{j=2}^{2N}\left|\phi_j(0)\right|^2= 0.
\end{equation}
We use in the following the basis introduced in the proof of Lemma \ref{lemV}.

First, since $\phi_{2j-1}(0)=0$ for $j=1,\cdots,N$, it is easy to see that there exists $a_{2j-1}$ such that
$$\phi_{2j-1}=a_{2j-1}\sin(z^2\cdot),\quad \forall j=1,\cdots,N.$$
Then, by the continuity at the interior nodes (\ref{p4}), we get
$$a_{2j-1}\sin(z^2 l_{2j-1})=0,\quad \forall j=1,\cdots,N.$$

Second, there exist $a_{2j}$, $b_{2j}$, $\tilde{a}_{2j}$ and $\tilde{b}_{2j}$ such that
$$\phi_{2j}=a_{2j}\sin(z\cdot)+b_{2j}\cos(z\cdot)+\tilde{a}_{2j}\sinh(z\cdot)+\tilde{b}_{2j}\cosh(z\cdot).
$$
By (\ref{p3}) and (\ref{sum=0}), we obtain
$$b_{2j}=\tilde{b}_{2j}=\tilde{a}_{2j}=0 \quad \hbox{and} \quad a_{2j}\sin(zl_{2j})=0 ,\quad \forall j=1,\cdots,N,
$$
since $z\neq0$.
Then, we have, with the notation introduced in (\ref{notcs}),
\begin{equation}\label{ajsj=0}
a_js_j=0,\quad j=1,\cdots, 2N.
\end{equation}

Moreover (\ref{p5}) gives
$$a_{2j}=\frac{1}{z}a_{2j-1}c_{2j-1},$$
and (\ref{p6}) yields
$$a_{2j+1}=za_{2j}c_{2j}.$$

By induction, we obtain, for all $j\geq2$,
\begin{equation}\label{inductionaj}a_j=z^{\epsilon_j}a_1c_{j-1}c_{j-2}\cdots c_1,\end{equation}
with $\epsilon_{2j}=-1$ and $\epsilon_{2j-1}=0$.
Therefore $a_1\neq0$ (otherwise $a_j=0$ for all j, and then $\phi=0$, which is impossible).
Now, by (\ref{ajsj=0}), we have $s_1=0$ and $c_1=\pm 1$. Then, since (\ref{inductionaj}) holds,
$a_2\neq0$ and $s_2=0$, again with (\ref{ajsj=0}). Then  $c_2=\pm 1$...
We see, by induction, that $s_j=0$ for all $j\in\left\{1,\cdots,2N\right\}$. 
Therefore, it suffices to have one $s_j\neq0$ for some $j\in\left\{1,\cdots,2N\right\}$ to obtain (\ref{sumneq0}).
It is the case if there exist $i,\,j\in\left\{1,\cdots,N\right\}$ such that
(\ref{condstab1}) or (\ref{condstab2}) hold.
\end{proof}

As a consequence of the previous lemma, we can prove the following proposition.
\begin{proposition}\label{propdecay0}
We have 
\begin{equation}\label{energy0}\displaystyle \lim_{t\rarrow + \infty}E(t)=0\end{equation}  
for all solution $\underline{u}$ of $(P_1)$ with
$(\underline{u}^0,\underline{u}^1)$ in ${\cal H}_1$ if and only if
(\ref{sumneq0}) holds for all eigenvectors $\underline{\phi}\in V_c$ of $(EP)$. 
Consequently, if there exist $i,\,j\in\left\{1,\cdots,N\right\}$ such that
(\ref{condstab1}) or (\ref{condstab2}) hold, then (\ref{energy0}) holds.
\end{proposition}

\begin{proof}
\framebox{$\Leftarrow$} Let us show that (\ref{sumneq0}) implies (\ref{energy0}). For that purpose we closely follow
\cite{tucsnak:03}.

First, we show that $\mathcal{A}_1$ has no eigenvalue on the imaginary axis.
If it is not the case, let $i\omega$ be an eigenvalue of $\mathcal{A}_1$ where $\omega\in\mathbb{R}^*$. Let $Z\in \mathcal{D}(\mathcal{A}_1)$ be an eigenvector associated with $i\omega$. Then $Z$ is of the form
$$Z=\left(\begin{array}{c}\underline{\phi}\\ i\omega \underline{\phi}\end{array}\right),$$
with 
\begin{equation}\begin{array}{c}
\partial_x^2\phi_{2j-1}=-\omega^2\phi_{2j-1},\quad j=1,\cdots,N,\\
\partial_x^4\phi_{2j}=\omega^2\phi_{2j},\quad j=1,\cdots,N.\label{equalityvp}\end{array}\end{equation} 
It is an immediate consequence of the identity $(i\omega I-\mathcal{A}_1)Z=0$.

We now take the inner product $\left\langle .,.\right\rangle_{\mathcal{H}_1}$ between $\mathcal{A}_1Z$ and $Z$.
By (\ref{dissipativeness}), we have
$$\Re\left(\left\langle\mathcal{A}_1Z,Z \right\rangle_{\mathcal{H}_1}\right)=-\omega^2\sum_{j=1}^{2N-1}\left|\phi_{j}(l_{j})\right|^2.$$
Since $Z$ is an eigenvector of $\mathcal{A}_1$ associated with $i\omega$ and $\omega\neq 0$, we obtain
$$\sum_{j=1}^{2N-1}\left|\phi_{j}(l_{j})\right|^2=0.
$$
Note that $Z$ satisfies the eigenvalue problem $(EP)$ and $Z$ belongs to $\mathcal{D}(\mathcal{A}_c)$,
since
$$P_{0,c}Z=\frac{1}{2i\pi}\oint_{\gamma}(\lambda I-\mathcal{A}_c)^{-1}Zd\lambda=\frac{1}{2i\pi}\oint_{\gamma}\frac{1}{\lambda-i\omega}Zd\lambda=0,
$$
(where we use $(\lambda I-\mathcal{A}_c)(\frac{1}{\lambda-i\omega}Z)=Z$ and where $\gamma$ is a simple closed curve enclosing only $0$),
and thus $Z=Z-P_{0,c}Z\in (I-P_{0,c})(V\times\prod_{j=1}^{2N}L^2(0,l_j))=\mathcal{H}_c$.
Then this contradicts (\ref{sumneq0}).
Therefore $\mathcal{A}_1$ has no eigenvalue on the imaginary axis.

Now, we can apply the main theorem of Arendt and Batty \cite{arendt:88}: Since $\sigma(\mathcal{A}_1)\cap i\mathbb{R}$ is
empty,
we obtain (\ref{energy0}).

\framebox{$\Rightarrow$} Let us show that (\ref{energy0}) implies (\ref{sumneq0}). 
For that purpose we use a contradiction argument.
Suppose that there exists an eigenvector  $\underline{\phi}\in V_c$ of $(EP)$ of associated eigenvalue $\lambda^2$ (where $\lambda=iz^2$, $z\in\mathbb{R}^{+*}$) such that
$$\sum_{j=1}^{2N-1}\left|\phi_{j}(l_{j})\right|^2=0.$$
Let us set
$$u(.,\,t)=\phi\cos(z^2 t).$$
Then $u$ is solution of $(P_1)$ and satisfies $$E(t)=E(0),$$
because
$$\phi_j(l_j)=0,\quad \forall j=1,\cdots, 2N.$$
This contradicts (\ref{energy0}).

\medskip
It suffices to use Lemma \ref{lemmsum} to finish the proof.
\end{proof}

Moreover, with the same method as previously, we are able to prove the decay to zero of the energy of solutions without restriction about the irrational  properties of the lengths.
\begin{proposition}
We have $\displaystyle \lim_{t\rarrow + \infty}E(t)=0$  for any solution of $(P_2)$ with
$(\underline{u}^0,\underline{u}^1)$ in ${\cal H}_2$.
\end{proposition}

\begin{proof}
As in the proof of Proposition \ref{propdecay0}, we can show that the energy of solutions of $(P_2)$ tends to zero if and only if 
\begin{equation}\label{sum2neq0}
\sum_{j=1}^{2N-1}\left|\phi_j(l_j)\right|^2+\sum_{j=1}^{N-1}\left(\left|\partial_x\phi_{2j}(l_{2j})\right|^2+\left|\partial_x\phi_{2j}(0)\right|^2\right)\neq0,
\end{equation}
for all eigenvectors $\underline{\phi}$ of $(EP)$. 
Let $\underline{\phi}$ be an eigenvector of $(EP)$ such that (\ref{sum2neq0}) is false.
By the same proof as Lemma \ref{lemmsum}, this implies that $\underline{\phi}=0$, which is impossible.
Then (\ref{sum2neq0}) holds and therefore the energy decays to 0.
\end{proof}

\begin{remark}
If we take the initial data in $V\times\prod_{j=1}^{2N}L^2(0,l_j)$, the energy of the solutions of $(P_1)$ and $(P_2)$ do not decay to $0$, since
$u=\phi$, where $(\phi,0)^t$ is an eigenvector of $\mathcal{A}_i$ ($i=1,\,2$) associated to the eigenvalue $0$, is solution of $(P_1)$ and $(P_2)$ with constant energy.
\end{remark}

%
%
%
%

\section{Stabilization result for $(P_2)$} \label{end}

We prove a decay result of the energy of system $(P_2)$, independently of the length of the strings and beams, for all regular initial data. 
In \cite{ammehjel}, the authors prove that the system described by $(P_2)$ is not exponentially stable in ${\cal H}_2$ with $N=1$ (i.e. with one string and one beam).
Therefore, in the general case (for $N\in\mathbb{N}^*$), we can not except to obtain an exponential decay for the energy of the solutions of $(P_2)$, but only a weaker decay rate, and in this general case, we prove a polynomial decay rate.
To obtain this, our technique is based on a frequency domain method and combines a contradiction argument with the multiplier technique to carry out a special analysis for the resolvent. 

The following theorem is a direct generalization of the result in \cite{ammehjel},
which we note, due to a mistake in the choice of
$\theta$, the decay rate in the following $\frac{ln^4(t)}{t^2}$ has
been written $\frac{ln^6(t)}{t^4}$ (corresponding to a choice of
$\theta = 1$ and not to $\theta = 1/2$).
\begin{theorem} \label{lr}
There exists a constant $C >0$ such that, for all $(\underline{u}^0,\underline{u}^1)\in{\cal D}({\cal A}_2)$, the solution of system $(P_2)$ satisfies the following estimate
\BEQ{EXPDECEXP3nb}
E(t)\le C \, \frac{\ln^{4}(t)}{t^2} \, \left\Vert (\underline{u}^0,\underline{u}^1) \right\Vert_{{\cal D}({\cal A}_2)}^2,
\FORALL t > 0.
\EEQ
\end{theorem}

 
\begin{proof}
We will employ the following frequency domain theorem for polynomial stability (see Liu-Rao \cite{rao}) of a $C_0$ semigroup
of contractions on a Hilbert space:

\begin{lemma}\label{lemrao}
A $C_0$ semigroup $e^{t{\cal L}}$ of contractions on a Hilbert space satisfies 
$$||e^{t{\cal L}}U_0|| \leq C \, \frac{\ln^{1 + \frac{1}{\theta}}(t)}{t^{\frac{1}{\theta}}} ||U_0||_{{\cal D}({\cal L})}$$
for some constant $C >0$ and for $\theta>0$ if 
\be 
\rho ({\cal L})\supset \bigr\{i \beta \bigm|\beta \in \rline \bigr\} \equiv i \rline, \label{1.8w} \ee 
and \be \limsup_{|\beta |\to \infty } \frac{1}{\beta^\theta} \, \|(i\beta -{\cal L})^{-1}\| <\infty, \label{1.9} 
\ee 
where $\rho({\cal L})$ denotes the resolvent set of the operator ${\cal L}$.
\end{lemma}

Then the proof of Theorem \ref{lr} is based on the following two lemmas.

\begin{lemma} \label{condsp}
The spectrum of ${\cal A}_2$ contains no point on the imaginary axis.
\end{lemma}

\begin{proof} Since ${\cal A}_2$ has compact resolvent, its spectrum $\s({\cal A}_2)$ only consists of eigenvalues of ${\cal A}_2$. We
will show that the equation 
\be {\cal A}_2 Z = i \b Z \label{1.10} \ee
with $Z= (\underline{y},\, \underline{v})^t \in \cD({\cal A}_2)$ and $\b\ne 0$ has only
the trivial solution.\\
\noindent By taking the inner product of (\ref{1.10}) with $Z$ and using 
\begin{multline}\label{1.7}
\Re
\left(<{\cal A}_2 Z,Z>_{{\cal H}_2} \right)
= - \sum_{j = 1}^N \left(\, \left| v_{2j}(0) \right|^2 + \left| \frac{d v_{2j}}{dx}(0) \right|^2 \right)\\ 
-\sum_{j = 1}^{N-1} \left(\left| v_{2j}(l_{2j}) \right|^2 + \left| \frac{d v_{2j}}{dx}(l_{2j}) \right|^2 \right), 
\end{multline}
we obtain that 
$$v_{2j}(0)=0,\, \frac{d v_{2j}}{dx}(0)=0, \, j=1,...,N \hbox{ and } v_{2j}(l_{2j})=0, \, \frac{d v_{2j}}{dx}(l_{2j})=0, \, j=1,...,N-1.$$ Next, we eliminate $\underline{v}$ in \rfb{1.10} to
get an ordinary differential equation: 
\be \left\{
\begin{array}{l}  
(\beta^2 y_{2j-1}+ \partial_x^2 y_{2j-1})(x)=0,\ x\in(0,l_{2j-1}), j = 1,...,N,\\ 
(\beta^2 y_{2j} - \partial_x^4 y_{2j})(x)=0,\ x \in (0,l_{2j}), j = 1,...,N,\\
y_1(0)=0,\ y_{2N}(l_{2N})=0,\ \partial_x^2 y_{2N}(l_{2N})=0,\\
\partial_x^2 y_{2j}(0)= 0, \, j = 1,...,N,\\
\partial_x^2 y_{2j}(l_{2j})= 0, \, j = 1,...,N-1,\\
y_j(l_j)=y_{j+1}(0),\,  j = 1,...,2N-1,\\
\partial_x^3 y_{2j}(0) + \partial_x y_{2j-1}(l_{2j-1})=0,\, j = 1,...,N,\\
\partial_x^3 y_{2j}(l_{2j}) + \partial_x y_{2j+1}(0)= 0,\, j = 1,...,N-1.
\end{array} \right. 
\label{1.11} 
\ee
Then, we can easily see that the only solution of the above system is the trivial one.
\end{proof}

The second lemma shows that (\ref{1.9}) holds with $\mathcal{L}=\mathcal{A}_2$ and $\theta=1$.

\begin{lemma}\label{lemresolvent}
The resolvent operator of $\mathcal{A}_2$ satisfies condition \rfb{1.9} for $\theta=1.$
\end{lemma}

\begin{proof}
Suppose that condition \rfb{1.9} is false with $\theta=1$. 
By the Banach-Steinhaus Theorem (see \cite{brezis}), there exists a sequence of real numbers $\beta_n \rightarrow +\infty$ and a sequence of vectors
$Z_n= (\underline{y}_{n},\, \underline{v}_{n} )^t\in {\cal D}({\cal A}_2)$ with $\|Z_n\|_{\cH_2} = 1$ such that 
\be 
|| \beta_n (i \b_n I - {\cal A}_2)Z_n||_{\cH_2} \rightarrow 0\;\;\;\; \mbox{as}\;\;\;n\rightarrow \infty, 
\label{1.12} \ee 
i.e., 
\be \beta_n^{1/2}\left(i \b_n y_{n} - v_{n}\right) \equiv f_{n}\rightarrow 0 \;\;\; \mbox{in}\;\; V, \label{1.13}\ee 
 \be
 \beta_n^{1/2} \left( i \b_n
v_{n,2j-1} - \frac{d^2 y_{n,2j-1}}{dx^2} \right) \equiv g_{n,2j-1} \rightarrow 0 \;\;\;
\mbox{in}\;\; L^2(0,l_{2j-1}),
\label{1.13b} \ee
 \be 
\beta_n^{1/2} \, \left( i \b_n v_{n,2j} + \frac{d^4 y_{n,2j}}{dx^4} \right) \equiv k_{n,2j} \rightarrow 0 \;\;\;
\mbox{in}\;\; L^2(0,l_{2j}), 
\label{1.14} \ee 
since $\beta_n^{1/2}\leq\beta_n$.

Our goal is to derive from \rfb{1.12} that $||Z_n||_{\cH_2}$ converges to zero, thus there is a contradiction. The proof is divided into four steps:

{\it First step.} We first notice that we have
\be 
|| \beta_n (i \b_n I - {\cal A}_2)Z_n||_{\cH_2} \ge |\Re \left(\langle \beta_n(i\beta_n I - {\cal A}_2)Z_n, Z_n\rangle_{\cH_2} \right)|. 
\label{1.15}
\ee
Then, by \rfb{1.7} and \rfb{1.12}, 
\be
\beta_n^{\frac{1}{2}}\,v_{n,2j}(0) \rightarrow 0, \quad  \beta_n^{\frac{1}{2}}\,\frac{d v_{n,2j}}{dx}(0) \rightarrow 0, \quad j=1,...,N
\label{1.16a}\ee
and 
\be
\beta_n^{\frac{1}{2}}\,v_{n,2j}(l_{2j}) \rightarrow 0, \quad \beta_n^{\frac{1}{2}}\,\frac{d v_{n,2j}}{dx}(l_{2j}) \rightarrow 0, \quad j=1,...,N-1.
\label{1.16} \ee 

This further leads, by (\ref{1.13}) and the trace theorem, to
\be
\left|\b_n \right|^{\frac{3}{2}}\, \left|y_{n,2j}(0) \right| \rightarrow 0, 
\quad |\beta_n|^{3/2} \, \left|\frac{dy_{n,2j}}{dx}(0) \right| \rightarrow 0, \quad j=1,...,N,
\label{1.17a} \ee
and
\be
\left|\b_n \right|^{\frac{3}{2}}\, \left|y_{n,2j}(l_{2j}) \right| \rightarrow 0, 
\quad |\beta_n|^{3/2} \, \left|\frac{dy_{n,2j}}{dx}(l_{2j}) \right| \rightarrow 0,  \quad j=1,...,N-1.
 \label{1.17} \ee
Moreover, since $Z_n\in{\cal D}({\cal A}_2)$ and thus satisfies (\ref{e6}), we have, by (\ref{1.16a}) and (\ref{1.16}),
\be
\left|\b_n \right|^{\frac{1}{2}}\,\left|\frac{d^2y_{n,2j}}{dx^2}(0)\right|\rightarrow 0, \,j=1,...,N,\quad \left|\b_n \right|^{\frac{1}{2}}\,\left|\frac{d^2y_{n,2j}}{dx^2}(l_{2j})\right|\rightarrow 0,\,j=1,...,N-1.
\label{4.50b}\ee
Then, note that, by continuity at the interior nodes and by (\ref{1.17a}) and (\ref{1.17}), we have
\be
\left|\b_n \right|^{\frac{3}{2}}\, \left|y_{n,2j-1}(0) \right| \rightarrow 0,\,j=2,...,N,\quad 
\left|\b_n \right|^{\frac{3}{2}}\, \left|y_{n,2j-1}(l_{2j-1}) \right| \rightarrow 0,\,j=1,...,N.
\label{4.50c}\ee

{\it Second step.} We now express $\underline{v}_n$ as a function of $\underline{y}_n$ from \rfb{1.13} and substitute it into
\rfb{1.13b}-\rfb{1.14} to get 
\be \beta_n^{1/2} \left( -\beta_n^2  y_{n,2j-1} - \frac{d^2
y_{n,2j-1}}{dx^2} \right) =  g_{n,2j-1} + i\beta_n f_{n,2j-1}, \quad j=1,...,N, \label{1.18} \ee
\be \beta_n^{1/2} \left( -\beta_n^2  y_{n,2j} + \frac{d^4
y_{n,2j}}{dx^4} \right) =  k_{n,2j} + i\beta_n f_{n,2j}, \quad j=1,...,N. 
\label{1.18b} 
\ee 

Next, we take the inner product of (\ref{1.18}) with $q_{2j - 1}(\cdot) \ds \frac{d y_{n,2j-1}}{dx}$
in $L^2(0,l_{2j-1})$ where $q_{2j - 1}\in C^1([0,l_{2j-1}])$ and $q_{2j - 1}(0)=0$. We obtain that
\begin{equation}\begin{array}{l}
\ds{\int_{0}^{l_{2j-1}} \beta_n^{1/2} \, \Bigr(-\b_n^2 y_{n,2j-1} - \frac{d^2 y_{n,2j-1}}{dx^2}
\Bigr) q_{2j-1}(x) \frac{d \bar{y}_{n,2j-1}}{dx} \, dx}\\
\ds{=\int_{0}^{l_{2j-1}} \Bigr( g_{n,2j-1} +i \b_n \, f_{n,2j-1} \Bigr) ~q_{2j-1}(x) \frac{d \bar{y}_{n,2j-1}}{dx} \, dx}\\
\ds{=\int_{0}^{l_{2j-1}} g_{n,2j-1} ~q_{2j-1}(x) \frac{d \bar{y}_{n,2j-1}}{dx} \, dx
 -i \int_{0}^{l_{2j-1}} q_{2j-1} \frac{d f_{n,2j-1}}{dx} ~ \beta_n \bar{y}_{n,2j-1} \, dx}\\
\ds{- i\int_{0}^{l_{2j-1}}  f_{n,2j-1} \frac{d q_{2j-1}}{dx} ~ \beta_n \bar{y}_{n,2j-1} \, dx
 + i f_{n,2j-1}(l_{2j-1}) q_{2j-1}(l_{2j-1}) \beta_n \bar{y}_{n,2j-1}(l_{2j-1}).}
\label{1.19}
\end{array}\end{equation}

It is clear that the right-hand side of (\ref{1.19}) converges to zero. Indeed, $f_{n,2j-1}$ and $g_{n,2j-1}$ converge to zero in $H^1(0,l_{2j-1})$ and $L^2(0,l_{2j-1})$ respectively, $\left\|Z_n\right\|_{\mathcal{H}_2}=1$ and (\ref{4.50c}) holds, and, finally, $\left|\beta_n y_{n,2j-1}\right|=\left|\frac{f_{n,2j-1}}{\beta_n^{1/2}}+v_{n,2j-1}\right|$ is bounded in $L^2(0,l_{2j-1})$.

By a straight-forward calculation,
\begin{multline*}
\Re \left\{\int_{0}^{l_{2j-1}} -\beta_n^2 y_{n,2j-1} \, q_{2j-1}
\frac{d \bar{y}_{n,2j-1}}{dx}\, dx \right\} = - \frac{1}{2} q_{2j-1}(l_{2j-1}) |\b_n y_{n,2j-1}(l_{2j-1})|^2 \\
+ \frac{1}{2} \int_{0}^{l_{2j-1}} \frac{dq_{2j-1}}{dx} |\b_ny_{n,2j-1}|^2dx 
\end{multline*}
and
\begin{multline*}
\Re \left\{ \int_{0}^{l_{2j-1}} - \,  \frac{d^2 y_{n,2j-1}}{dx^2}\, q_{2j-1}
\frac{d \bar{y}_{n,2j-1}}{dx}
\, dx \right\}  = - \frac{1}{2} \, q_{2j-1}(l_{2j-1})\left|\frac{d y_{n,2j-1}}{dx}(l_{2j-1})\right|^2 \\
+\frac{1}{2} \int_{0}^{l_{2j-1}} \left| \frac{d y_{n,2j-1}}{dx} \right|^2 
\, \frac{d q_{2j-1}}{dx} dx.
\end{multline*}

We then take the real part of (\ref{1.19}), and (\ref{4.50c}) leads to
\begin{multline}
 \int_{0}^{l_{2j-1}} \frac{d q_{2j-1}}{dx}  \left|\beta_n y_{n,2j-1} \right|^2dx +
\int_{0}^{l_{2j-1}}   \frac{d q_{2j-1}}{dx} \left|\frac{d y_{n,2j-1}}{dx}
\right|^2 dx\\ 
- q_{2j-1}(l_{2j-1}) \left|\frac{d y_{n,2j-1}}{dx}(l_{2j-1})\right|^2 \rightarrow 0. 
\label{1.20}
\end{multline}

Similarly, we take the inner product of (\ref{1.18b}) with $q_{2j} (\cdot) \ds \frac{d y_{n,2j}}{dx}$ in $L^2(0,l_{2j})$ with $q_{2j} \in
C^3([0,l_{2j}])$ and $q_{2j}(l_{2j})=0$. We then repeat the above procedure. Since 
\begin{multline*}
\int_{0}^{l_{2j}} \left|\frac{dy_{n,2j}}{dx} \right|^2 dx = - \frac{1}{i \beta_n} \int_{0}^{l_{2j}} v_{n,2j} \frac{d^2 \bar{y}_{n,2j}}{dx^2} - \frac{1}{i \beta_n} \int_{0}^{l_{2j}} (i\beta_n \, y_{n,2j} - v_{n,2j}) \frac{d^2\bar{y}_{n,2j}}{dx^2} dx \\
- \frac{d\bar{y}_{n,2j}}{dx}(0)\, y_{n,2j}(0)
+\frac{d\bar{y}_{n,2j}}{dx}(l_{2j})\, y_{n,2j}(l_{2j}),
\end{multline*}
then, from the boundedness of $v_{n,2j}$, $i\beta_n y_{n,2j} - v_{n,2j}$, $\frac{d^2 y_{n,2j}}{dx^2}$ in $L^2(0,l_{2j})$ and (\ref{1.17a})-(\ref{1.17}), $\frac{dy_{n,2j}}{dx}$ converges to zero in $L^2(0,l_{2j})$. 
This will give, after some calculations,
\begin{multline}
\int_{0}^{l_{2j}}  \frac{d q_{2j}}{dx} |\beta_ny_{n,2j}|^2 dx +
 \int_{0}^{l_{2j}} 3\, \frac{d q_{2j}}{dx}  \left|
\frac{d^2 y_{n,2j}}{dx^2} \right|^2 dx \\
- 2 \Re\left( \frac{d^3 y_{n,2j}}{dx^3}(0) q_{2j}(0) \frac{d\bar{y}_{n,2j}}{dx}(0)\right)\rightarrow 0. 
\label{1.21}
\end{multline}

{\it Third step.} Next, we show that $\ds \frac{d y_{n,2j-1}}{dx}(l_{2j-1})$ and $\ds \frac{d^3 y_{n,2j}}{dx^3}(0)$ converge to zero. 
We take the inner product of \rfb{1.18b} with $\frac{1}{\beta_n^{1/2}}e^{-\beta_n^{1/2}x}$ in $L^2(0,l_{2j})$.
We have, with \rfb{1.18b},
\be\begin{array}{rcl}
\ds{ \int_{0}^{l_{2j}}\left(-\beta_n^2y_{n,2j}+\frac{d^4u_{n,2j}}{dx^4}\right)e^{-\beta_n^{1/2}x}dx}
&=&\ds{\int_{0}^{l_{2j}}\frac{1}{\beta_n^{1/2}}k_{n,2j}e^{-\beta_n^{1/2}x}dx}\\
&&\ds{+i\int_{0}^{l_{2j}}\beta_n^{1/2}f_{n,2j}e^{-\beta_n^{1/2}x}dx}.
\end{array}\label{tend0}\ee
It is clear that the first term of the right hand side of (\ref{tend0}) tends to zero by (\ref{1.14}). 
Moreover, by integration by parts,
$$\int_{0}^{l_{2j}}\beta_n^{1/2}f_{n,2j}e^{-\beta_n^{1/2}x}dx
=\int_{0}^{l_{2j}}\frac{d f_{n,2j}}{dx} e^{-\beta_n^{1/2}x}dx-f_{n,2j}(l_{2j})e^{-\beta_n^{1/2}l_{2j}}+f_{n,2j}(0),
$$
which tends to zero since $f_{n,2j}$ tends to zero in $H^2$ and by the trace theorem. 

This leads to 
\be  \int_{0}^{l_{2j}}\left(
  \beta_n^{2} e^{-\beta_n^{1/2} x}\, y_{n,2j} -  e^{-\beta_n^{1/2} x}\,\frac{d^4 y_{n,2j}}{dx^4}\, \right) dx \rightarrow 0. \label{1.22bb}
\ee

Performing four integrations by parts in the second term on the left-hand side of \rfb{1.22bb}, we obtain
\begin{multline}
\int_{0}^{l_{2j}}\left(\beta_n^{2} e^{-\beta_n^{1/2} x}\, y_{n,2j} -  e^{-\beta_n^{1/2} x}\,\frac{d^4 y_{n,2j}}{dx^4}\, \right) dx 
= \frac{d^3 y_{n,2j}}{dx^3}(0)+ \beta_n^{1/2}\,\frac{d^2 y_{n,2j}}{dx^2}(0) \\
+\beta_n\,\frac{d y_{n,2j}}{dx}(0) + \beta_n^{3/2} y_{n,2j}(0) + o(1),
\label{1.22}\end{multline}
with (\ref{1.17})-(\ref{4.50b}) and since 
$$\begin{array}{rcl}
\ds{\left|\frac{d^3 y_{n,2j}}{dx^3}(l_{2j})e^{-\beta_n^{1/2}l_{2j}}\right|^2}
&\leq& \ds{e^{-2\beta_n^{1/2}l_{2j}}\int_{0}^{l_{2j}}\left|\frac{d^4 y_{n,2j}}{dx^4}(x)\right|^2dx}\\
&\leq& \ds{e^{-2\beta_n^{1/2}l_{2j}}\int_{0}^{l_{2j}}\left|\frac{k_{n,2j}}{\beta_n^{1/2}}-i\beta_n v_{n,2j}\right|^2dx}\\
&\leq&\ds{\frac{2}{\beta_n}e^{-2\beta_n^{1/2}l_{2j}}\int_{0}^{l_{2j}}\left|k_{n,2j}\right|^2dx}\\
&&\ds{+ 2\beta_n^2 e^{-2\beta_n^{1/2}l_{2j}}\int_{0}^{l_{2j}}\left|v_{n,2j} \right|^2dx\rightarrow 0, }
\end{array}$$
because $\left\|Z_n\right\|_{\mathcal{H}_2}=1$.

Thus, according to \rfb{1.17a} and (\ref{4.50b}), we simplify (\ref{1.22}) to 
\be
\label{1.24vb}
\frac{d^3 y_{n,2j}}{dx^3}(0) \rightarrow 0.
\ee
Consequently, since $Z_n\in\mathcal{D}(\mathcal{A}_2)$ and thus satisfies (\ref{e4}), we obtain
\be
\label{1.24vbb}
\frac{d y_{n,2j-1}}{dx}(l_{2j-1}) \rightarrow 0.
\ee

Then, (\ref{1.17a}) and (\ref{1.24vb}) lead to
\be
\frac{d \bar{y}_{n,2j}}{dx}(0) \, \frac{d^3 y_{n,2j}}{dx^3}(0) \rightarrow 0.
\label{1.26}
\ee

In view of \rfb{1.24vbb}-\rfb{1.26}, we simplify \rfb{1.20} and \rfb{1.21} to
\begin{equation}
 \int_{0}^{l_{2j-1}}  \frac{d q_{2j-1}}{dx}|\beta_ny_{n,2j-1}|^2dx +
 \int_{0}^{l_{2j-1}}   \frac{d q_{2j-1}}{dx} \left|\frac{d y_{n,2j-1}}{dx} \right|^2 dx
 \rightarrow 0,  \label{1.28}
 \ee
 \be
 \int_{0}^{l_{2j}}  \frac{d q_{2j}}{dx} |\beta_ny_{n,2j}|^2dx +
 \int_{0}^{l_{2j}} 3\, \frac{d q_{2j}}{dx}  \left|
\frac{d^2 y_{n,2j}}{dx^2} \right|^2 dx
\rightarrow 0 \label{1.21b}\ee 
respectively.

{\it Fourth step.} Finally, we choose $q_{2j-1}$ and $q_{2j}$ such that $\ds \frac{d q_{2j-1}}{dx}$ is strictly positive and $\ds \frac{d q_{2j}}{dx}$ is strictly negative.
This can be done by taking $$ q_{2j-1}(x)= e^{x}-1, \qquad q_{2j}(x) = e^{(l_{2j} -x)}-1. $$ Therefore, (\ref{1.28}) and
(\ref{1.21b}) imply 
\be \|\b_ny_{n,2j-1}\|_{L^2(0,l_{2j-1})}\rightarrow0, \, \|\b_ny_{n,2j}\|_{L^2(0,l_{2j})}\rightarrow
0,\, \|(y_{n,2j-1},y_{n,2j})_{j\in\left\{1,\cdots,N\right\}}\|_{V}\rightarrow 0. 
\label{1.31}
\ee 

In view of (\ref{1.13}), we also get 
\be \|v_{n,2j-1}\|_{L^2(0,l_{2j-1})}\rightarrow 0, \quad \|v_{n,2j}\|_{L^2(0,l_{2j})} \rightarrow 0,
\label{1.32}\ee 
which clearly contradicts $\left\|Z_n\right\|_{{\cal H}_2}=1$. 
\end{proof}

The two hypothesis of Lemma \ref{lemrao} are proved by Lemma \ref{condsp} and Lemma \ref{lemresolvent}.
Then (\ref{EXPDECEXP3nb}) holds.
The proof of Theorem \ref{lr} is then finished.
\end{proof}


\begin{thebibliography}{99}
\bibitem{ammari1} K. Ammari and M. Tucsnak, Stabilization of Bernoulli-Euler
beams by means of a pointwise feedback force,
{\em SIAM J. Control. Optim.,} {\bf 39} (2000), 1160-1181. 
\bibitem{ammari2} K. Ammari, A. Henrot and M. Tucsnak, Asymptotic behaviour of the solutions and optimal location of the actuator for the
pointwise stabilization of a string, {\em Asymptotic Analysis,} {\bf 28} (2001),
215-240. 
\bibitem{ammari3} K. Ammari, Z. Liu and M. Tucsnak, Decay rates for a beam with pointwise force and moment feedback, {\em Mathematics of
Control, Signals, and systems,} {\bf 15} (2002), 229-255. 
\bibitem{ammari4} K. Ammari and M. Jellouli, Remark in stabilization of tree-shaped networks of
strings, {\em Appl. Maths.,} {\bf 4} (2007), 327-343. 
\bibitem{ammari5} K. Ammari, Asymptotic behaviour of some elastic planar networks of Bernoulli-Euler beams, {\em Appl. Anal.,} {\bf 86} (2007), 1529-1548. 
\bibitem{ammari} K. Ammari and M. Tucsnak, Stabilization of second order evolution equations by a class of unbounded feedbacks, {\em ESAIM Control Optim. Calc. Var.,} {\em ESAIM Control Optim. Calc. Var,} {\bf 6} (2001), 361-386. 
\bibitem{amjel} K. Ammari and M. Jellouli, Stabilization of star-shaped networks of strings, {\em Diff. Integral. Equations,}
{\bf 17} (2004), 1395-1410. 
\bibitem{amjellk} K. Ammari, M. Jellouli and M. Khenissi, Stabilization of generic trees of strings, 
{\em J. Dyn. Cont. Syst.,} {\bf 11} (2005), 177-193. 
\bibitem{ammehjel} K. Ammari, M. Jellouli and M. Mehrenberger, Feedback stabilization of a coupled string-beam system, {\em Netw. Heterog. Media.}, {\bf 4} (2009), 19-34.
\bibitem{arendt:88} W.~Arendt and C.~J.~K. Batty, Tauberian theorems and stability of one-parameter semigroups,
{\em Trans. Amer. Math. Soc.}, {\bf 305}(1988), 837-852.
\bibitem{banks2} H. T. Banks, R. C. Smith and Y. Wang, {\em Smart Materials Structures}, Wiley, 1996.
\bibitem{brezis} H. Brezis, {\em Analyse Fonctionnelle, Th\'eorie et Applications},
Masson, Paris, 1983. 
\bibitem{jvb} J. von Below, Classical solvability of linear parabolic equations on networks, {\em J. Diff. Eq.}, \bf 72 \rm (1988), 316-337.
\bibitem{dagerzuazua} R.~D{\'a}ger and E.~Zuazua, {\em Wave propagation, observation and control in {$1\text{-}d$} flexible multi-structures}, volume 50 of Math\'ematiques \& Applications (Berlin), Springer-Verlag, 2006.
\bibitem{Kato} T. Kato, {\em Perturbation theory for linear operators}, Reprint of the 1980 Edition, Springer-Verlag, Berlin, 1995.
\bibitem{lagnese} J. Lagnese, G. Leugering and E. J. P. G. Schmidt, {\em Modeling, Analysis of dynamic elastic multi-link structures}, Birkh\"auser, Boston-Basel-Berlin, 1994. 
\bibitem{rao} Z. Liu and B. Rao, Characterization of polynomial decay rate for the solution of linear evolution equation, {\em Z. Angew. Math. Phys.,} {\bf 56} (2005), 630-644. 
\bibitem{mer} D. Mercier, Spectrum analysis of a serially connected Euler-Bernouilli beams  problem,
{\em Netw. Heterog. Media.}, {\bf 4} (2009), 709-730.
\bibitem{merreg1} D. Mercier, V. Régnier, Spectrum of a network of Euler-Bernoulli beams, {\em J. Math. Anal. and Appl.}, {\bf 337} (2007), 174-196.
\bibitem{merreg2} D. Mercier, V. Régnier, Control of a network of Euler-Bernoulli beams, {\em J. Math. Anal. and Appl.}, {\bf 342} (2008), 874-894.
\bibitem{merreg3} D. Mercier, V. Régnier, Boundary controllability of a chain of serially connected Euler-Bernoulli beams with interior masses, 
{\em Collect. Math}, {\bf 60}  (2009), 307-334.
\bibitem{valein:07} S.~Nicaise and J.~Valein, Stabilization of the wave equation on 1-{D} networks with a delay term in the nodal feedbacks, {\em Netw. Heterog. Media}, {\bf 2} (2007), 425--479.
\bibitem{paulsen} W.H. Paulsen, The exterior matrix method for sequentially coupled fourth-order equations, {\em J. of Sound and Vibration}, 
{\bf 308} (2007), 132-163.
\bibitem{Pazy} A. Pazy, {\em Semigroups of linear operators and applications
to partial differential equations}, Springer, New York, 1983.
\bibitem{tucsnak:03} M.~Tucsnak and G.~Weiss, How to get a conservative well-posed linear system out of thin air.
  {II}. {C}ontrollability and stability, {\em SIAM J. Control Optim.}, 42(3):907--935, 2003.
\bibitem{tucsnakbook} M.~Tucsnak and G.~Weiss, Observation and control for operator semigroups, {\em Birkh\"auser Advanced Texts: Basler Lehrb\"ucher}, Birkh\"auser Verlag, Basel, 2009.
\bibitem{xumasto} G.Q. Xu and N.E. Mastorakis, Stability of a star shaped coupled networks of strings and beams,
{\em WSEAS, Proceeding of the 10th WSEAS international conference on Technique and Computations}, Technical University of Sofia (Bulgaria), 2008.
\bibitem{zhangxumasto} K.T. Zhang, G.Q. Xu and N.E. Mastorakis, Stability of a complex network of {E}uler-{B}ernoulli beams,
{\em WSEAS Trans. Syst.}, {\bf 8} (2009), 379--389.
\end{thebibliography}
\end{document}